
\def\title#1{{\titlefont\noindent #1\bigskip}}

\def\author#1{{\largefont\noindent #1}\medskip}

\def\beginlinemode{\endmode
 \begingroup\obeylines\def\endmode{\par\endgroup}}
\let\endmode=\par

\newbox\theaddress
\def\address{\smallskip\beginlinemode\parindent 0in\getaddress}
{\obeylines
\gdef\getaddress #1 
 #2
 {#1\gdef\addressee{#2}%
   \global\setbox\theaddress=\vbox\bgroup\raggedright%
    \everypar{\hangindent2em}#2
   \def\endaddress{\egroup\endgroup \copy\theaddress \medskip}}}

\def\thanks#1{\footnote{}{\eightpoint #1}}

\long\def\Abstract#1{{\eightpoint\narrower\vskip\baselineskip\noindent
#1\smallskip}}

\def\skipfirstword#1 {}

\def\ir#1{\csname #1\endcsname}

\newdimen\currentht
\newbox\droppedletter
\newdimen\droppedletterwdth
\newdimen\drophtinpts
\newdimen\dropindent

\def\irrnSection#1#2{\edef\tttempcs{\ir{#2}}
\currentht-\pagetotal\advance\currentht by-\ht\footins
\advance\currentht by\vsize
\ifdim\currentht<1.5in\par\vfill\eject\else\vbox to.25in{\vfil}\fi
{\largefont\noindent{}\expandafter\skipfirstword\tttempcs. #1}
\vskip\baselineskip }

\def\irSubsection#1#2{\edef\tttempcs{\ir{#2}}
\vskip\baselineskip\penalty-3000
{\bf\noindent \expandafter\skipfirstword\tttempcs. #1}
\vskip\baselineskip}

\def\irSubsubsection#1#2{\edef\tttempcs{\ir{#2}}
\vskip\baselineskip\penalty-3000
{\bf\noindent \expandafter\skipfirstword\tttempcs. #1}
\vskip\baselineskip}

\def\References{\vbox to.25in{\vfil}\noindent{}{\bf References}
\vskip\baselineskip\par}

\def\baselinebreak{\par \ifdim\lastskip<\baselineskip
         \removelastskip\penalty-200\vskip\baselineskip\fi}

\long\def\prclm#1#2#3{\baselinebreak
\noindent{\bf \csname #2\endcsname}:\enspace{\sl #3\par}\baselinebreak}

\def\Prf{\noindent{\bf Proof}: }

\def\rem#1#2{\baselinebreak\noindent{\bf \csname #2\endcsname}:\enspace }

\def\qed{{\hfill$\diamondsuit$}\vskip\baselineskip}

\def\bibitem#1{\par\indent\llap{\rlap{\bf [#1]}\indent}\indent\hangindent
2\parindent\ignorespaces}

\long\def\eatit#1{}

\def\leftheadlinetext{}
\def\rightheadlinetext{}

\def\leftheadline{{\eightrm\folio\hfil \leftheadlinetext\hfil}}
\def\rightheadline{{\eightrm\hfil\rightheadlinetext\hfil\folio}}

\headline={\ifnum\pageno=1\hfil\else
\ifodd\pageno\rightheadline\else\leftheadline\fi\fi}

\def\tenpoint{\def\rm{\fam0\tenrm}
\textfont0=\tenrm \scriptfont0=\sevenrm \scriptscriptfont0=\fiverm
\textfont1=\teni \scriptfont1=\seveni \scriptscriptfont1=\fivei
\def\mit{\fam1} \def\oldstyle{\fam1\teni}
\textfont2=\tensy \scriptfont2=\sevensy \scriptscriptfont2=\fivesy
\def\cal{\fam2}
\textfont3=\tenex \scriptfont3=\tenex \scriptscriptfont3=\tenex
\def\it{\fam\itfam\tenit} 
\textfont\itfam=\tenit
\def\sl{\fam\slfam\tensl} 
\textfont\slfam=\tensl
\def\bf{\fam\bffam\tenbf} 
\textfont\bffam=\tenbf \scriptfont\bffam=\sevenbf
\scriptscriptfont\bffam=\fivebf
\def\tt{\fam\ttfam\tentt} 
\textfont\ttfam=\tentt
\normalbaselineskip=12pt
\setbox\strutbox=\hbox{\vrule height8.5pt depth3.5pt  width0pt}%
\normalbaselines\rm}

\def\eightpoint{\def\rm{\fam0\eightrm}%
\textfont0=\eightrm \scriptfont0=\sixrm \scriptscriptfont0=\fiverm
\textfont1=\eighti \scriptfont1=\sixi \scriptscriptfont1=\fivei
\def\mit{\fam1} \def\oldstyle{\fam1\eighti}%
\textfont2=\eightsy \scriptfont2=\sixsy \scriptscriptfont2=\fivesy
\def\cal{\fam2}%
\textfont3=\tenex \scriptfont3=\tenex \scriptscriptfont3=\tenex
\def\it{\fam\itfam\eightit} 
\textfont\itfam=\eightit
\def\sl{\fam\slfam\eightsl} 
\textfont\slfam=\eightsl
\def\bf{\fam\bffam\eightbf} 
\textfont\bffam=\eightbf \scriptfont\bffam=\sixbf
\scriptscriptfont\bffam=\fivebf
\def\tt{\fam\ttfam\eighttt} 
\textfont\ttfam=\eighttt
\normalbaselineskip=9pt%
\setbox\strutbox=\hbox{\vrule height7pt depth2pt  width0pt}%
\normalbaselines\rm}

\def\largefont{\def\rm{\fam0\largerm}
\textfont0=\largerm \scriptfont0=\largescriptrm \scriptscriptfont0=\tenrm
\textfont1=\largei \scriptfont1=\largescripti \scriptscriptfont1=\teni
\def\mit{\fam1} \def\oldstyle{\fam1\teni}
\textfont2=\largesy 
\def\cal{\fam2}
\def\it{\fam\itfam\largeit} 
\textfont\itfam=\largeit
\def\sl{\fam\slfam\largesl} 
\textfont\slfam=\largesl
\def\bf{\fam\bffam\largebf} 
\textfont\bffam=\largebf 
\def\tt{\fam\ttfam\largett} 
\textfont\ttfam=\largett
\normalbaselineskip=17.28pt
\setbox\strutbox=\hbox{\vrule height12.25pt depth5pt  width0pt}%
\normalbaselines\rm}

\def\titlefont{\def\rm{\fam0\titlerm}
\textfont0=\titlerm \scriptfont0=\largescriptrm \scriptscriptfont0=\tenrm
\textfont1=\titlei \scriptfont1=\largescripti \scriptscriptfont1=\teni
\def\mit{\fam1} \def\oldstyle{\fam1\teni}
\textfont2=\titlesy 
\def\cal{\fam2}
\def\it{\fam\itfam\titleit} 
\textfont\itfam=\titleit
\def\sl{\fam\slfam\titlesl} 
\textfont\slfam=\titlesl
\def\bf{\fam\bffam\titlebf} 
\textfont\bffam=\titlebf 
\def\tt{\fam\ttfam\titlett} 
\textfont\ttfam=\titlett
\normalbaselineskip=24.8832pt
\setbox\strutbox=\hbox{\vrule height12.25pt depth5pt  width0pt}%
\normalbaselines\rm}

\nopagenumbers

\font\eightrm=cmr8
\font\eighti=cmmi8
\font\eightsy=cmsy8
\font\eightbf=cmbx8
\font\eighttt=cmtt8
\font\eightit=cmti8
\font\eightsl=cmsl8
\font\sixrm=cmr6
\font\sixi=cmmi6
\font\sixsy=cmsy6
\font\sixbf=cmbx6

\font\largerm=cmr12 at 17.28pt
\font\largei=cmmi12 at 17.28pt
\font\largescriptrm=cmr12 at 14.4pt
\font\largescripti=cmmi12 at 14.4pt
\font\largesy=cmsy10 at 17.28pt
\font\largebf=cmbx12 at 17.28pt
\font\largett=cmtt12 at 17.28pt
\font\largeit=cmti12 at 17.28pt
\font\largesl=cmsl12 at 17.28pt

\font\titlerm=cmr12 at 24.8832pt
\font\titlei=cmmi12 at 24.8832pt
\font\titlesy=cmsy10 at 24.8832pt
\font\titlebf=cmbx12 at 24.8832pt
\font\titlett=cmtt12 at 24.8832pt
\font\titleit=cmti12 at 24.8832pt
\font\titlesl=cmsl12 at 24.8832pt

\tenpoint



\def\manyby{\hbox to.75in{\hrulefill}}
\hsize 6.5in 
\vsize 9.2in

\tolerance 3000
\hbadness 1500
\vbadness 1500

\def\item#1{\par\indent\indent\llap{\rlap{#1}\indent}\hangindent
2\parindent\ignorespaces}

\def\itemitem#1{\par\indent\indent
\indent\llap{\rlap{#1}\indent}\hangindent
3\parindent\ignorespaces}

\def\trans{H1}
\def\vanc{H2}
\def\bowdoin{H3}
\def\ravello{H4}
\def\mtnwest{H5}
\def\antican{H6}
\def\fatpts{H7}

\def\binom#1#2{\hbox{$\left(\matrix{#1\cr #2\cr}\right)$}}
\def\C#1{\hbox{$\cal #1$}}
\def\s{\hbox{$\scriptstyle\cal S$}}
\def\r{\hbox{$\scriptstyle\cal R$}}

\def\pr#1{\hbox{{\bf P}${}^{#1}$}}
\def\leftheadlinetext{Brian Harbourne}
\def\rightheadlinetext{Fat Point Ideal Generation}

\title{The Ideal Generation Problem for Fat Points}

\author{Brian Harbourne}

\address
Department of Mathematics and Statistics
University of Nebraska-Lincoln
Lincoln, NE 68588-0323
email: bharbourne@unl.edu
WEB: http://www.math.unl.edu/$\sim$bharbour/
\smallskip
March 26, 1997\endaddress
\vskip-\baselineskip

\thanks{\vskip -6pt
\noindent This work benefitted from a National Science Foundation grant.
\smallskip
\noindent 1991 {\it Mathematics Subject Classification. } 
Primary 13P10, 14C99. 
Secondary 13D02, 13H15.
\smallskip
\noindent {\it Key words and phrases. }  Ideal generation 
conjecture, symbolic powers, syzygies,
resolution, fat points, maximal rank, plane, Cremona.\smallskip}

\vskip\baselineskip
\Abstract{Abstract: This paper is concerned with the problem of 
determining up
to graded isomorphism the modules in a minimal free resolution
of a fat point subscheme $Z=m_1p_1+\cdots+m_rp_r\subset\pr2$
for general points $p_1,\ldots,p_r$.}
\vskip\baselineskip

\irrnSection{Introduction}{intro}
We always work over an arbitrary algebraically closed field $k$.
This paper is concerned with determining
the number $\nu_t(I(Z))$ of elements 
in each degree $t$ of any minimal set of homogeneous generators in
the ideal $I(Z)\subset k[\pr2]$ defining a fat point 
subscheme $Z=m_1p_1+\cdots+m_rp_r\subset\pr2$, where
$p_1,\ldots,p_r\in\pr2$ are general. Given the Hilbert
function of $I(Z)$, this is equivalent up to graded isomorphism to
determining the modules in a minimal free resolution of $I(Z)$. 

As discussed further below, the above problem
has been solved for subschemes $Z=p_1+\cdots+p_r\subset\pr2$,
for general points $p_i$.
The solution rests on showing in such cases that $I(Z)$
has the maximal rank property: given a graded ideal $I$ in
a polynomial ring $R$ graded in the usual way by degree,
we say that $I$ has the {\it maximal rank property\/}
if the multiplication maps $\mu_t(I):I_t\otimes R_1\to I_{t+1}$
have maximal rank (i.e., are injective or surjective)
for every $t$ (where the subscript $t$ denotes homogeneous
components of degree $t$).

Since fat point subschemes commonly fail to have
the maximal rank property, it has been unclear
what sort of answer to the general problem can be expected.
In this paper we suggest an asymptotic solution.
In particular, fixing points $p_1,\ldots,p_r\in\pr2$,
we define an equivalence relation, 
Cremona equivalence, on fat point
subschemes $Z=m_1p_1+\cdots+m_rp_r$, and, 
if the points are strongly nonspecial (see \ir{bnds})
and $Z$ is expectedly good (a property defined below
giving control over the Hilbert function of $I(Z)$, and which
holds in all known cases for general points $p_i$),
we show in \ir{applythm} for all but finitely many subschemes $Z$ in 
each Cremona equivalence class that $I(Z)$ satisfies the maximal rank property.

There nevertheless remains the problem of understanding failures of 
the maximal rank property. Sometimes this is easy.
For example, given a graded ideal $I\subset R$,
define $\alpha(I)$ to be
the least degree among nonzero homogeneous elements of $I$,
and define $\beta(I)$ to be the least degree $t$ such that 
the elements of $I_t$ have no nontrivial common divisor.
It is easy to see that $\mu_t$ cannot be injective for
$t>\alpha$ and cannot be surjective for $t=\beta-1$,
so having $\alpha<\beta-1$ guarantees that $\mu_{\beta-1}$
fails to have maximal rank and thus
that $I$ does not have the maximal rank property.
On the other hand, failures of $\mu_\beta$ to have maximal rank
are more mysterious,
and, in fact, by \ir{cokfact} and \ir{myomega} 
the general problem of determining numbers of generators
for expectedly good fat point subschemes reduces
to determining the rank of $\mu_\beta$.

For an expectedly good fat point subscheme $Z$
with $\alpha<\beta$, Fitchett [Fi] 
shows that the greatest common 
divisor of $I(Z)_\alpha$ determines the rank of $\mu_\beta$. This work
gives a geometric explanation for the possible failure of maximal rank
of $\mu_\beta$ in the case that $\beta>\alpha$,
in addition to determining bounds on the rank of $\mu_\beta$.

What is still lacking is 
a general understanding of why $\mu_\beta$ could
fail to have maximal rank in the case that $\alpha=\beta$.
Our result \ir{applythm} on Cremona equivalence actually shows
in the expectedly good case not only that the
maximal rank property holds asymptotically but
that $\alpha=\beta$ holds asymptotically as well.
However, $\mu_\beta$ can fail to have maximal rank
even if $\alpha=\beta$,
and we study this phenomenon
in the case of {\it uniform\/} fat point subschemes
(i.e., subschemes $Z=m(p_1+\cdots+p_r)$). 
For example, from our results in 
\ir{smallr} it follows that:

\prclm{Corollary}{mxrkfails}{Let $p_1,\ldots,p_r$
be $r\le 9$ general points of \pr2 and let 
$I=I(m(p_1+\cdots+p_r))$. Then
$\alpha(I)=\beta(I)$ but $\mu_\beta(I)$ fails to have maximal
rank if and only if: $r=7$, $m=3l$ and $3\le l\le 7$;
or $r=8$, $m=6l$ and $9\le l\le 16$; or
$r=8$, $m=6l+1$ and $6\le l\le 13$.}

We give the proof in \ir{corproof}.
Our results of \ir{smallr} also explicitly
determine the modules in a  minimal free resolution
of $I(m(p_1+\cdots+p_r))$ for any $m$ and for 
any $r\le 9$ general points of \pr2.

For $r>9$ the question remains open,
but in \ir{VIGC} we propose for uniform subschemes
that the failures in \ir{mxrkfails} above are the 
only failures for any $r$ general points. We also provide some
evidence for this in \ir{bnds}, 
using Campanella-like bounds (viz. \ir{mycamp}, cf. [Cam])
to verify a number of cases of the conjecture
for expectedly good fat point subschemes.

We will use the following notational convention.
A divisor on a surface $X$ will be denoted with the typeface {\tt C}.
Its class in the divisor class group $\hbox{Cl}(X)$ (of divisors
modulo linear equivalence) will be
denoted $C$, and the corresponding line bundle in $\hbox{Pic}(X)$
will be \C C. In certain special cases, we will also use lower case 
letters to denote divisor classes, and $\C O_X(F)$ to denote the
line bundle corresponding to a class $F$. Finally, in certain
instances it will be convenient not to discriminate
between a divisor class and its corresponding line bundle,
which we may do, for example, by writing $H^i(X, F)$ in place of 
the strictly correct $H^i(X, \C O_X(F))$.

\irSubsection{Previous Work}{prevwork}
To put the results of this paper into the context of other
recent work, let $I\subset R$ be an ideal (where
$R=k[x_0,\ldots,x_n]$ is a polynomial ring),
homogeneous with respect to the usual grading (in which
each indeterminate $x_i$ has degree 1 and constants have degree 0).

A typical approach to understanding $I$ begins with 
its Hilbert function (which gives the $k$-vector 
space dimension $\hbox{dim }I_t$
of each graded component $I_t$ as a function of the degree $t$).
Next one looks at the number $\nu_t(I)$ of 
elements of degree $t$ in any minimal set of homogeneous generators;
this gives the first module in a minimal free resolution for $I$.
Finally, one considers the successive syzygy modules in a minimal
free resolution.

In trying to elucidate principles governing the behavior of these
aspects of ideals of $R$, it is natural to regard $R$ as the homogeneous
coordinate ring of the projective space 
\pr n of dimension $n$, and to begin
with ideals associated to subvarieties or subschemes of \pr n.
(The reader will recall the usual bijection $X\mapsto I(X)$
from closed subschemes of \pr n to saturated homogeneous ideals of $R$.)

Points being the geometrically simplest subschemes, one is naturally 
attracted to studying ideals of the form $I(m_1p_1+\cdots+m_rp_r)$,
for distinct points $p_1,\ldots,p_r\in \pr n$ 
and nonnegative integers $m_i$, not all 0,
where $I(m_1p_1+\cdots+m_rp_r)$ denotes the 
homogeneous ideal generated by all forms
which vanish at each point $p_i$ with multiplicity at least $m_i$.
Following Geramita, the corresponding 
subscheme $m_1p_1+\cdots+m_rp_r$ is called
a {\it fat point subscheme} and its ideal 
$I(m_1p_1+\cdots+m_rp_r)$ is called
a {\it fat point ideal}.

For general points $p_1,\ldots,p_r$, the ideals $I(p_1+\cdots+p_r)$ 
have been studied extensively (viz.,
[HS], [Lor] and [EP]).
In this situation, the Hilbert function is 
known trivially (each point imposing 
independent conditions on forms of each 
degree until no forms of that degree remain)
so attention has focused on numbers of generators and on resolutions.
Of particular interest here is the 
{\it Ideal Generation Conjecture\/} (IGC)
of [GO] and [GGR]:

\prclm{Ideal Generation Conjecture}{igc}{The ideal $I(Z)$ 
has the maximal rank property 
for any general set $Z=p_1+\cdots+p_r$ of $r$ points in \pr n.}

To see its relevance, note 
for any homogeneous ideal $J\subset R$ that 
$\nu_{t+1}(J)$ is the dimension of the cokernel of the multiplication
map $\mu_t(J): J_t\otimes R_1\to J_{t+1}$ defined for $f\in J_t$ by 
$f\otimes x_i\mapsto x_if$. If the Hilbert function of $J$ is known
(and thus the dimensions of $J_t\otimes R_1$ and $J_{t+1}$), then
the rank of $\mu_t(J)$ determines $\hbox{dim cok }\mu_t(J) = \nu_{t+1}(J)$.

Although this conjecture remains open in general,
it has been verified in various cases
(see [Bl], [GM], [HS], [HSV], [Lor], [O], [Ra]), 
including $n=2$ for all $r$ [GGR]. In addition,
on \pr2 a minimal free resolution of $I=I(Z)$
is of the form $0\to F_1\to F_0\to I\to 0$, where
$F_0=\oplus_t R[-t]^{\nu_t(I)}$. Thus given the number
$\nu_t(I)$ of generators for each $t$ and the Hilbert function
of $I$, one knows the Hilbert function of $F_1$ and hence
one knows $F_1$ itself. In particular, the problem of determining
the minimal free resolution of $I(Z)$ on \pr2 reduces to
determining the Hilbert function and numbers $\nu_t(I)$
of generators, and is thus completely solved for
any general set $Z=p_1+\cdots+p_r\subset\pr2$.

Much less is known or even conjectured
in the situation $m_1p_1+\cdots+m_rp_r$ of fat points,
in which the coefficients $m_i$ need not be at most 1.
Most work either restricts $r$, $n$ or the 
coefficients $m_i$. For example,
[Cat] completely works out
the minimal free resolution for any $m_i$ for $r<6$
general points and $n=2$
([Fi] extends this to $r=6$), while [A], 
[AH1], [AH2], [AH3], [Hi], [Ch] 
determine the Hilbert function 
for any $r$ and $n$ if each $m_i$ is at most 2 
and [CM] for any $r$ with $n=2$ and $m_i$ small and 
nearly constant. 
Some steps toward understanding 
the Hilbert function of generally situated fat points in \pr n
have been taken (viz. [I]), but only for \pr 2 has a conjecture for
the Hilbert function of any generally situated finite set of 
fat points been suggested
(first in [\vanc] and later equivalent 
variants in [Hi], [Gi] and [\ravello]).

\irSubsection{\pr2 and its Blowings up}{blups}
Thus only for \pr2 do we have a putative Hilbert function for generally
situated fat points, and this begs the questions of what we should expect
for the numbers of generators (and hence for the 
minimal free resolution), given the expected behavior
for Hilbert functions. Although it is an open question whether
the expected behavior is always obtained, it can in many cases be
verified.

We now discuss this in more detail.
To do so, we must consider
surfaces obtained by blowing up points of \pr2. 
In particular, let $p_1,\ldots,p_r$ be distinct points of \pr2. 
Let $\pi:X\to \pr2$ 
be the morphism obtained by blowing up 
each point $p_i$. Let ${\tt E}_i$ denote the exceptional
divisor of the blow up of $p_i$, and let $e_i$ denote its divisor class.
Let $e_0$ denote the pullback to $X$
of the class of a line in \pr2; the 
classes $e_0,\ldots,e_r$ comprise a ${\bf Z}$-basis
of $\hbox{Cl}(X)$. Note that this basis, which we 
call an {\it exceptional configuration}, 
is completely determined by $\pi$ and in turn
determines $\pi$. Also, recall that 
$\hbox{Cl}(X)$ supports an intersection form with respect to which
the basis $e_0,\ldots,e_r$ is orthogonal, satisfying
$-1=-e_0^2=e_1^2=\cdots=e_r^2$, and that the canonical class $K_X$ 
of $X$ is $K_X=-3e_0+e_1+\cdots+e_r$.
Recall that a divisor class is {\it numerically effective\/} if its
intersection with every effective divisor is nonnegative,
and that a prime divisor {\tt C} on $X$
with $C^2=-1=C\cdot K_X$ is smooth and rational, called a
{\it $(-1)$-curve}, or an {\it exceptional curve}. We refer to its class 
$C$ as a {\it $(-1)$-class} or an {\it exceptional class}. 
It is known precisely which classes
are exceptional classes, when $p_1,\ldots,p_r$ are sufficiently general.

To establish the connection to fat points, consider  
a fat point subscheme $Z=m_1p_1+\cdots+m_rp_r\subset\pr2$. 
Let $X$ be obtained by blowing up 
each point $p_i$ and let $e_0,\ldots,e_r$ be the corresponding 
exceptional configuration. Let $F_d$ denote the class
$de_0-m_1e_1-\cdots-m_re_r$. Since $e_0$ corresponds to the pullback
$\pi^*(\C O_{\pr2}(1))$ of the class of a line, 
we have for each $d$ and $i$ a natural isomorphism 
of $H^i(X,\C F_d)$ with 
$H^i(\pr2,\pi_*(\C O_X(-m_1e_1-\cdots-m_re_r))\otimes
\C O_{\pr2}(d))=H^i(\pr2,\C I_Z(d))$.
In particular, the homogeneous coordinate ring 
$R=\bigoplus_{d\ge 0}H^0(\pr2,\C O_{\pr2}(d))$ can be identified
with $\bigoplus_{d\ge 0}H^0(X, de_0)$, and the
homogeneous ideal 
$I(Z)=\bigoplus_{d\ge 0}H^0(\pr2,\C I_Z(d))$ in $R$  
can be identified with 
$\bigoplus_{d\ge 0}H^0(X,\C F_d)$. Moreover, 
under these identifications, the 
natural homomorphisms $H^0(X,\C F_d)\otimes 
H^0(X, e_0)\to H^0(X,\C F_{d+1})$
and $I(Z)_d\otimes R_1\to I(Z)_{d+1}$ correspond, so
the dimension $\nu_{d+1}$ of the cokernel of the latter is equal to
the dimension of the cokernel of the former. 

Now, suppose $F=F_d$ is the class of an effective 
divisor. By taking $N$ to comprise 
the components of negative self-intersection in the fixed locus of $|F|$,
we can write $F=H+N$, where $H$ and $N$
are the classes of effective divisors, $H$ is numerically 
effective with $h^0(X,\C F)=h^0(X,\C H)$, and $N$ is a sum of
prime divisors of negative self-intersection with
$h^0(X,\C N)=1$. If the points $p_1,\ldots,p_r$ are general,
in all known cases it is true that $h^1(X,\C H)=0$ and that
$N$ is a sum of multiples of classes of disjoint exceptional curves
disjoint from a general element of $|H|$. In such a case,
since the exceptional classes are known, we can explicitly determine
$N=-\sum (E\cdot F)E$ (where the sum is over all
exceptional classes $E$ with $E\cdot F<0$), 
and hence the value $h^0(X, \C F)=
(H^2-H\cdot K_X)/2+1$ of the Hilbert function of $I(Z)$
in degree $d$. Assuming the foregoing behavior always holds, 
we can also explicitly determine
whether $F_d$ is the class of an effective divisor 
(see [\trans], [\mtnwest]).
The point of this paper is to assume the foregoing situation holds,
and study the consequences for determining numbers of generators.
Toward this end, we make the following definition.

\rem{Definition}{expgood} Let
$Z=m_1p_1+\cdots+m_rp_r\subset\pr2$ 
be a fat point subscheme, let $X$ be the blowing
up of the points $p_i$ and let 
$F_t=te_0-m_1e_1-\cdots-m_re_r$. Then we say $Z$
is {\it expectedly good\/}
if $F_{\alpha(I(Z))}=H+N$, where $H$ is numerically effective
and $N$ is a 
nonnegative sum of exceptional classes with
$h^0(X,\C O_X(F_{\alpha(I(Z))}))=h^0(X, \C H)$,
$h^1(X,\C H)=0$ and $h^0(X, \C N)=1$. 
(It easily follows that $H\cdot N=0$ and thus that
$N=-\sum (E\cdot F)E$, where the sum is over all
exceptional classes $E$ with $E\cdot F<0$.)
We also say that
the points $p_1,\ldots,p_r\subset\pr2$ are
{\it expectedly good\/} if the only prime divisors on $X$
of negative self-intersection are exceptional curves
and if for every effective and numerically effective
divisor ${\tt C}$ we have $h^1(X,\C O_X({\tt C}))=0$.
\vskip\baselineskip

Note that if $p_1,\ldots,p_r$ are expectedly good, then so is any
$Z=m_1p_1+\cdots+m_rp_r$, and, if $Z$ is expectedly good,
one only needs to know the classes of exceptional curves and 
the coefficients $m_i$ in order to compute
the Hilbert function of $I(Z)$. 

By [\mtnwest], $r\le 8$ general points $p_1,\ldots,p_r\in \pr 2$ 
are expectedly good, and each $Z=m_1p_1+\cdots+m_9p_9$
is expectedly good for general points $p_1,\ldots,p_9$.
Any 9 sufficiently general points, by which we mean the complement
of a countable union of closed conditions (which is nonempty
unless $k$ is the algebraic closure of a finite field), are also expectedly good.
On the other hand, three or more collinear
points, or six or more on a conic, or the nine base points of a 
cubic pencil are not expectedly good.
Whether 10 or more sufficiently general points are expectedly 
good is unknown, but they are expected to be, and conjectures
to this effect have been put forward (viz. [\vanc], [Hi], [Gi] and [\ravello]).
Moreover, many specific examples of expectedly good fat point subschemes
$Z=m_1p_1+\cdots+m_rp_r$ are known with $r>9$. 

\irSubsection{A Generalized IGC}{genIGC}
Let us say that the Uniform Maximal Rank Property (UMRP) on \pr n holds
for $r$ if, for each $m>0$, the 
maximal rank property for $I(mp_1+\cdots+mp_r)$ holds
for general points 
$p_1,\ldots,p_r$ of \pr n. Let us also say that
the Restricted Uniform Maximal Rank Property (RUMRP) on \pr n holds
for $r$ if $\mu_{\beta(I(mp_1+\cdots+mp_r))}$ 
has maximal rank for each $m>0$
for general points $p_1,\ldots,p_r$ 
of \pr n. We will show
in \ir{smallr} that:

\prclm{Theorem}{data}{Let $r\le 9$. Then the UMRP on \pr 2
holds if and only if $r$ is 1, 4, or 9, and the RUMRP holds
if and only if $r$ is not 7 or 8.}

For general points $p_1,\ldots,p_r\in\pr2$, failures of maximal rank 
seem for uniform $Z$ to be confined to small $r$. 
For example, the failure of the UMRP on \pr2 when $r$ is 
a nonsquare less than 9 is, by \ir{abfail},
guaranteed by the existence of
uniform abnormal curves for such $r$.
(Following Nagata [N1], a curve ${\tt C}\subset\pr 2$ of degree $d$
whose multiplicity at each point 
$p_i$ is at least $m_i$ is said to be {\it abnormal\/} if 
$d\sqrt{r}<m_1+\cdots+m_r$, and {\it uniform\/}
if $m_1=\cdots=m_r$.) But Nagata [N1] proves that
no abnormal curves occur for $r$ generic points when $r$ is a square, 
and he [N2] conjectures
that none occur for $r>9$. Additional evidence that
we present in \ir{bnds} also suggests that 
the RUMRP may hold on \pr2 for $r>9$. Moreover, for $r>9$ 
expectedly good points, RUMRP implies UMRP by \ir{mycampcor}.
This prompts us, with some temerity perhaps, to 
propose a generalized IGC for uniform fat points:

\prclm{Conjecture}{VIGC}{The UMRP on \pr2 holds for all $r>9$.}

This also suggests the following question:

\prclm{Question}{ques}{Is there an $N$ depending on $n$, 
such that the UMRP holds on \pr n for each $r\ge N$?}

\irrnSection{Background on Surfaces}{surf}
For the rest of this paper,
$R$ will denote the homogeneous coordinate ring
$R=k[x_0,x_1,x_2]$ of \pr2 (over any algebraically closed field  $k$).
Let $X$ be obtained by blowing up distinct
points $p_1,\dots,p_r\in\pr2$ and let $e_0,\ldots,e_r$ be the corresponding 
exceptional configuration. Let $F_d$ denote the class
$de_0-m_1e_1-\cdots-m_re_r$ and let $Z=m_1p_1+\cdots+m_rp_r$. 

Following [Mu], we will
denote the kernel of $H^0(X,\C F_d)\otimes 
H^0(X, e_0)\to H^0(X,\C F_{d+1})$
by $\C R(\C F_d,e_0)$ and the cokernel by $\C S(\C F_d,e_0)$;
it is then convenient to denote their dimensions
by $\r(\C F_d,e_0)$ and $\s(\C F_d,e_0)$.
Note that to say that $I(Z)_d\otimes R_1\to I(Z)_{d+1}$,
or equivalently that $H^0(X,\C F_d)\otimes 
H^0(X,e_0)\to H^0(X,\C F_{d+1})$,
has maximal rank is just to say that $[\r(\C F_d,e_0)][\s(\C F_d,e_0)]=0$.

First we have:

\prclm{Proposition}{Mumford}{Let $T$ be a closed 
subscheme of projective space,
let \C A and \C B be coherent sheaves on 
$T$ and let \C C be the class of an
effective divisor ${\tt C}$ on $T$.
\item{(a)} If the restriction homomorphisms 
$H^0(T, \C A)\to H^0({\tt C},\C A\otimes\C O_{\hbox{\eighttt C}})$
and $H^0(T, \C A\otimes\C B)\to H^0(C,\C A\otimes
\C B\otimes\C O_{\hbox{\eighttt C}})$ are surjective
(for example, if $h^1(T,\C A\otimes\C C^{-1})=0=
h^1(T,\C A\otimes\C C^{-1}\otimes\C B)$),
then we have an exact sequence
$$\eqalign{0\to & \C R(\C A\otimes\C C^{-1},\C B)\to
\C R(\C A,\C B)\to\C R(\C A\otimes\C O_{\hbox{\eighttt C}},\C B)\to\cr
&\C S(\C A\otimes\C C^{-1},\C B)\to\C S(\C A,\C B)\to
\C S(\C A\otimes\C O_{\hbox{\eighttt C}},\C B)\to0.\cr}$$
\item{(b)} If $H^0(T, \C B)\to H^0(C,\C B\otimes
\C O_{\hbox{\eighttt C}})$ is surjective
(for example, if $h^1(T,\C B\otimes\C C^{-1})=0$), then 
$\C S(\C A\otimes\C O_{\hbox{\eighttt C}},\C B)=\C S(\C A\otimes
\C O_{\hbox{\eighttt C}},\C B\otimes\C O_{\hbox{\eighttt C}})$.
\item{(c)} If $T$ is a smooth curve of genus 
$g$, and \C A and \C B are line 
bundles of degrees at least $2g+1$ and $2g$, respectively, then 
$\C S(\C A,\C B)=0$.}

\Prf See [Mu] for (a) and (c); we leave (b) as 
an easy exercise for the reader.\qed

Let $F$ be the class of an effective divisor ${\tt F}$ on a surface $X$. 
Let $F=H+N$ be a Zariski decomposition of $F$ (i.e.,
$h^0(X,\C F)=h^0(X,\C H)$ and $h^0(X,\C N)=1$; for example,
$N$ could be the class of the fixed part of the linear 
system $|{\tt F}|$ and
then $H=F-N$ would be the free part of $F$). 
The following lemma allows us to 
reduce a consideration of $\C S(\C F,e_0)$
to one of $\C S(\C H,e_0)$. 

\prclm{Lemma}{cokfact}{Let $e_0,\ldots,e_r$ 
be the exceptional configuration
corresponding to a blowing up $\pi:X\to\pr2$ at distinct points
$p_1,\ldots,p_r$, and let $F$ be a divisor class on $X$. 
If $F$ is not the class of an effective divisor, then 
$\s(\C F,e_0)=h^0(X, F+e_0)$.
If $F$ is the class of
an effective divisor, let $F=H+N$ be a Zariski decomposition; 
then $\s(\C F,e_0)=[h^0(X,F+e_0)-h^0(X,H+e_0)]+\s(\C H,e_0)$.}

\Prf See Lemma 2.10 of [\fatpts]. \qed

\rem{Remark}{recall} To determine $\nu_t(I(Z))$ for each $t$ for
some fat point subscheme $Z=m_1p_1+\cdots+m_rp_r$ of \pr2,
by \ir{cokfact} it is enough on the 
blow up $X$ of \pr2 at $p_1,\ldots,p_r$ to
determine $h^0(X, de_0-m_1e_1-\cdots-m_re_r)$ 
for every $d$, and, for each $d$ such that
$h^0(X, de_0-m_1e_1-\cdots-m_re_r)>0$, to 
determine: the free part $H$ of 
$de_0-m_1e_1-\cdots-m_re_r$; $\s(\C H,e_0)$; and $h^0(X, H+e_0)$.
(It is not hard to see that being able to compute
$h^0(X, \C F)$ for any $F$ also lets one find the free part of any $F$
whenever $F$ is the class of an effective divisor. And
once one knows $\nu_t(I(Z))$ for all $t$, one also knows the modules
in a minimal free resolution 
$0\to F_1\to F_0\to I(Z)\to 0$ of $I(Z)$; viz. \ir{example}.) 

In the case of any $r\le 9$ points, the 
results of [\antican] provide a solution
to determining $h^0(X,\C F)$ for any class $F$, and thus
to finding the free part of $F$ 
when $h^0(X,\C F)>0$. For $r\le 9$ general points, these results 
are well known and can, in any case, 
be recovered from [\antican] or [\trans]; for the 
reader's convenience we recall the facts relevant to a uniform class
$F$ in the case of $r$ general points 
of \pr2, first for $r\le 8$, and then for $r=9$.
(A class $F$ on a blowing up $X$ of \pr2 at distinct points
$p_1,\ldots,p_r$ will be called a {\it uniform\/} class if
$F=de_0-m(e_1+\cdots+e_r)$ for some nonnegative integers $d$ and $m$.)

Let $X$ be the blowing up of $r\le 8$ 
general points of \pr2. If $F$ is uniform and
if it is the class of an effective divisor, then 
the fixed part $N$ is also uniform,
equal to $-\sum (E\cdot F)E$, where the 
sum is over all classes $E$ of $(-1)$-curves 
with $E\cdot F<0$.
The classes of the $(-1)$-curves are known; 
up to permutation of the indices, 
they are (see Section 26 of [Ma]):
$e_1$, $e_0-e_1-e_2$, $2e_0-(e_1+\cdots+e_5)$, 
$3e_0-(2e_1+e_2+\cdots+e_7)$, 
$4e_0-(2e_1+2e_2+2e_3+e_4+\cdots+e_8)$, 
$5e_0-(2e_1+\cdots+2e_6+e_7+e_8)$, 
and $6e_0-(3e_1+2e_2+\cdots+2e_8)$. Now 
one can show that $N=0$ if $r=1$ or 4;
otherwise, $N$ is
a nonnegative multiple of: $e_0-e_1-e_2$ if 
$r=2$; $3e_0-2e_1-2e_2-2e_3$
if $r=3$; $2e_0-(e_1+\cdots+e_5)$ for $r=5$; 
$12e_0-5(e_1+\cdots+e_6)$, $r=6$;
$21e_0-8(e_1+\cdots+e_7)$, $r=7$; or 
$48e_0-17(e_1+\cdots+e_8)$, $r=8$.
It also follows that a uniform class 
$de_0-m(e_1+\cdots+e_r)$ is the class of an effective
divisor if and only if $d\ge \epsilon_r m$, 
where $\epsilon_1=\epsilon_2=1$,
$\epsilon_3=3/2$, $\epsilon_4=\epsilon_5=2$, 
$\epsilon_6=12/5$, $\epsilon_7=21/8$,
and $\epsilon_8=48/17$.

Now, the free part of the class of an effective divisor 
is numerically effective, and, if $X$ is any
blowing up of \pr2 at 8 or fewer points, 
general or not, then  ([\mtnwest], [\antican])
any numerically effective class
$F$ on $X$ is the class of an effective 
divisor with no fixed components 
and has $h^1(X,\C F)=h^2(X,\C F)=0$, hence 
$h^0(X,\C F)=(F^2-K_X\cdot F)/2+1$
by Riemann--Roch for surfaces.

Finally, let $r=9$. Nine general points 
of \pr2 always lie on a smooth cubic
curve, so more generally let $X$ be the blowing up of 
any $r=9$ distinct points 
of a smooth cubic curve ${\tt C}'$ on \pr2. 
Then $-K_X=3e_0-e_1-\cdots-e_9$ is numerically effective, the class 
of a smooth elliptic curve ${\tt C}$, the proper transform of ${\tt C}'$.
If $F$ is uniform, we can write $F=te_0-sK_X$ for uniquely determined
integers $t$ and $s$, with $s\ge 0$.
Moreover, $F$ is the class of an effective divisor
if and only if $t$ is also nonnegative, in which case $h^2(X,\C F)=0$,
hence $h^0(X,\C F)=(F^2-K_X\cdot F)/2+1+h^1(X,\C F)$. 
In addition, if $t>0$, then
$F$ is fixed part free and $h^1(X,\C F)=0$.
If, however, $t=0$, things are more delicate.
If the restriction of $\C O_X(-K_X)$ to ${\tt C}$ has infinite order
in $\hbox{Pic}({\tt C})$, let $a=0$.
Otherwise, let $l$ be the order of the 
restriction of $\C O_X(-K_X)$ to ${\tt C}$,
and define $a$ and $b$ via $s=al+b$ where $0\le b<l$. Then 
$h^0(X,-sK_X)=a+1$ and $h^1(X,-sK_X)=a$.
(Note for an algebraically closed field 
$k$ which is not the algebraic closure of a finite field, that
for sufficiently general---i.e., a nonempty complement of
a countable union of closed conditions---points 
$p_1,\ldots,p_9$ no nonzero power of $\C O_X(-K_X)$
restricts trivially to ${\tt C}$. If $k$ is the algebraic 
closure of a finite field, however, then 
the restriction of $\C O_X(-K_X)$ to ${\tt C}$ always has finite order.)
\vskip\baselineskip

The next result will be helpful in verifying failure of the UMRP.
Call a uniform class $E=de_0-m(e_1+\cdots+e_r)$ 
on a blowing up $X$ of \pr2 at distinct points
$p_1,\ldots,p_r$ {\it abnormal\/} if $E$ is the class of 
an effective divisor with $d<\sqrt{r}m$
(note that this is equivalent to $E^2<0$).

\prclm{Proposition}{abfail}{Let $X$ be a blowing up
of $r$ distinct points $p_1,\ldots,p_r$ of \pr2. If $X$ has
a uniform abnormal class $E$, then for some positive integers
$n$ and $m$, $I(m(p_1+\cdots+p_r))_n
\otimes R_1\to I(m(p_1+\cdots+p_r))_{n+1}$
does not have maximal rank.}

\Prf Since $E$ is the class of an 
effective divisor of negative self-intersection,
we can find positive integers $a$ and 
$b$ such that $ae_0+bE$ has nontrivial fixed part
but such that $(a+1)e_0+bE$ has trivial fixed part. Now,
$ae_0+bE=ne_0-m(e_1+\cdots+e_r)$ for some positive $n$ and $m$.
Since $a>0$, $H^0(X,ae_0)\otimes 
H^0(X,e_0)\to H^0(X,(a+1)e_0)$ is not injective,
hence neither is $H^0(X,ae_0+bE)\otimes 
H^0(X,e_0)\to H^0(X,(a+1)e_0+bE)$.
Since $(a+1)e_0+bE$ is fixed component free but $ae_0+bE$ is not,
we see $H^0(X,ae_0+bE)\otimes H^0(X,e_0)\to H^0(X,(a+1)e_0+bE)$
is also not surjective. Thus 
$H^0(X,ae_0+bE)\otimes H^0(X,e_0)\to H^0(X,(a+1)e_0+bE)$,
and hence $I(m(p_1+\cdots+p_r))_n
\otimes R_1\to I(m(p_1+\cdots+p_r))_{n+1}$, do not have maximal rank.\qed

The following result is well known (see 
Proposition 3.7 of [DGM]) and follows easily
by appropriately applying \ir{Mumford} (or by noting that 
$\tau_Z+1$ is just the regularity $\sigma(I(Z))$ of $I(Z)$).

\prclm{Lemma}{myomega}{Let $e_0,\ldots,e_r$ be the classes
corresponding to a blowing up $X\to\pr2$ at
distinct points $p_1,\ldots, p_r$. Let $Z=m_1p_1+\cdots+m_rp_r$,
and let $F_d$ denote $de_0-m_1e_1-\cdots-m_re_r$.
If $\omega_Z$ is the degree of a generator of 
greatest degree in a minimal
set of homogeneous generators of $I(Z)$ 
(equivalently, $\mu_d$ is surjective
for $d\ge \omega_Z$ but not for $d=\omega_Z-1$) 
and if $\tau_Z$ is the least integer
$t$ such that $h^1(X,\C F_t)=0$, then  $\omega_Z\le \tau_Z+1$.
In particular, $\C S(\C F_t,e_0)=0$ for $t>\tau_Z$.}

We now determine Campanella-like bounds (cf. [Cam]). 
Let $F=a_0e_0-a_1e_1-\cdots-a_re_r$,
with $a_i>0$ for all $i$, 
be the class of an effective divisor
on a blow up $X$ of \pr2 at distinct points $p_1,
\ldots,p_r\in\pr2$.
Let $h=h^0(X,\C F)$, $l_i=h^0(X,F-(e_0-e_i))$, 
and $q_i=h^0(X,F-e_i)$. 

\prclm{Lemma}{mycamp}{Given the multiplication map
$\mu: H^0(X,\C F)\otimes H^0(X, e_0)\to H^0(X, F + e_0)$
and $0<i\le r$, we have:
$$\hbox{max}(l_i,3h-h^0(X,F+e_0)) \le \hbox{dim ker $\mu$} 
\le l_i+q_i.$$}

\Prf For specificity, take $i=1$.
Let $x$ ($y$ and $z$, resp.) be the equation of 
the line through $p_2$ and $p_3$ (resp.,
$p_1$ and $p_3$, and $p_1$ and $p_2$). Let $L$ be the image of
$\Gamma(e_0-e_1)$ in $\Gamma(e_0)$, where $\Gamma$ is 
the global sections functor. Thus $L$ can be regarded
as the vector space span of $y$ and $z$, making 
$z\Gamma(\C F)+y\Gamma(\C F)$ the image
of $\Gamma(\C F)\otimes L$ under $\mu$. It has dimension $2h-l_1$
since $z\Gamma(\C F)\cap y\Gamma(\C F)=zy\Gamma(F-(e_0-e_1))$,
where we regard the intersection as 
taking place in $\Gamma(((F\cdot e_0)+1)e_0)$.
Therefore, $l_1 \le \hbox{dim ker $\mu$}$. But since
$\Gamma(\C F)\otimes\Gamma(e_0)$ has dimension $3h$ and $\mu$ maps into
$H^0(X, F + e_0)$, it is clear that we also have 
$3h-h^0(X,F+e_0) \le \hbox{dim ker $\mu$}$
and hence $\hbox{max}(l_1,3h-h^0(X,F+e_0)) \le \hbox{dim ker $\mu$}$.

To bound \hbox{dim ker $\mu$} above, note that all elements of 
$z\Gamma(\C F)+y\Gamma(\C F)$ correspond to forms
on \pr2 that vanish at $p_1$ to order at least $a_1+1$.
Thus $(y\Gamma(\C F)+z\Gamma(\C F))\cap x\Gamma(\C F)$
lies in the image of $x\Gamma(\C F-e_1)$ under the natural
inclusion $x\Gamma(F-e_1)\subset x\Gamma(\C F)$, so 
$\hbox{dim Im $\mu$}\ge (2h-l_1) + (h-q_1)$
hence $\hbox{dim ker $\mu$}\le l_1+q_1$.
\qed

\prclm{Corollary}{Firstcampcor}{Let $F$ and $\mu$ be as in \ir{mycamp},
let $d=F\cdot e_0$ and assume $h^1(X,F)=0$.
\item{(a)} Then $\mu$ has maximal rank if 
and only if $\hbox{max}(0,2h-d-2)=\hbox{dim ker $\mu$}$.
\item{(b)} Moreover, $\hbox{max}(0,2h-d-2)
\le \hbox{dim ker $\mu$} \le l_1+q_1$.
\item{(c)} If $h^1(X,F-(e_0-e_1))=h^1(X,F-e_1)=0$, 
then $l_1+q_1=2h-d-2$.}

\Prf We use the notation of \ir{mycamp}.

(a) Clearly, $\mu$ has maximal rank if 
and only if $\hbox{max}(0,3h-h^0(X,F+e_0))=
\hbox{dim ker $\mu$}$. But $h^1(X, F)=0$ 
(and hence $h^1(X, F+e_0)=0$),
so by Riemann--Roch we compute $h^0(X,F+e_0)
=h+d+2$. Thus $3h^0(X,F)-h^0(X,F+e_0)=
2h-d-2$ and the result follows.

(b) This follows by the proof of (a) and 
by \ir{mycamp}.

(c) Let $m=F\cdot e_1$. 
Since $h^1(X,F-e_1)=0$, taking ${\tt E}$ to be the effective 
divisor whose class is $e_1$, 
the exact sheaf sequence $0\to \C O_X(F-e_1)
\to \C F\to \C O_{\hbox{\eighttt E}}\otimes \C F\to 0$
is exact on global sections, so $h=h^0(X, F)=h^0(X,F-e_1)+
h^0({\tt E},\C O_X(F)\otimes\C O_{\hbox{\eighttt E}})=q_1+m+1$.

Since $h^1(X,F-(e_0-e_1))=0$, taking ${\tt C}$ to be a
general effective 
divisor whose class is $e_0-e_1$, 
the exact sheaf sequence $0\to \C O_X(F-(e_0-e_1))
\to \C F\to \C O_{\hbox{\eighttt C}}\otimes \C F\to 0$
is exact on global sections. Computing dimensions we find
$h=l_1+d+1-m$ so $2h-d-2=l_1+(h-m-1)=l_1+q_1$. 
\qed

\irrnSection{Applying the Bounds}{bnds}
Let $X$ be the blow up of \pr2 at distinct points
$p_1,\ldots,p_r\in\pr2$. Let $e_0,\ldots,e_r$ be the corresponding 
exceptional configuration, and 
define the {\it roots\/} $\rho_0=e_0-e_1-e_2-e_3$,
$\rho_i=e_i-e_{i+1}$, $i>0$. Reflections $s_i(x)=x+(x\cdot \rho_i)\rho_i$
through each $\rho_i$ define 
intersection form-preserving involutions of $\hbox{Cl}(X)$,
generating a subgroup $W$ (infinite for $r>8$), 
called the {\it Weyl\/} group, of the orthogonal group on
$\hbox{Cl}(X)$. Let us say 
that $p_1,\ldots,p_r$ are {\it strongly nonspecial\/}
if $h^0(X,\C F)=h^0(X,\C O_X(wF))$ 
for all $w\in W$ and $F\in \hbox{Cl}(X)$.
This is somewhat stronger than but implies Nagata's condition of
being nonspecial for Cremona transformations [N1].
And just as points which are independent generic points
over the prime field are nonspecial for Cremona
transformations [N1], they are also strongly nonspecial.
Nor is it hard to check that expectedly good points 
are strongly nonspecial. As a further example,
over any algebraically closed ground field $k$, sufficiently
general smooth points of a cuspidal cubic ${\tt C}'$ are strongly nonspecial.
(By sufficiently general, taking $X$ to be the blow up of
\pr2 at the points and ${\tt C}$ to be the proper 
transform to $X$ of ${\tt C}'$,
we mean such that the kernel of the induced homomorphism
$\hbox{Pic}(X)\to \hbox{Pic}({\tt C})$ 
is trivial in characteristic 0 or is $pK_X^\perp$ in characteristic
$p$, where $K_X^\perp$ is the subgroup of elements
\C F with $F\cdot K_X=0$. For justification, see Example 3.4 of
[\bowdoin], and use [\trans].)

We will obtain some asymptotic results that essentially say
that some property holds for all but finitely many elements
of a Weyl group orbit. The next lemma determines some properties 
of these orbits, including that they tend to be infinite.

\prclm{Lemma}{Worbits}{Let $F\ne 0$ be a numerically effective class
on the blowing up $X$ of
strongly nonspecial points $p_1,\ldots,p_r\in\pr2$, where $e_0,\ldots,
e_r$ is the corresponding exceptional configuration. 
\item{(a)} The orbit $WF$ under the Weyl group action is infinite
if and only if $r>9$, or $r=9$ but $F\ne -lK_X$ for any $l>0$.
\item{(b)} The class $wF-e_0$ is the class of an effective divisor
for at most finitely many elements $wF$ of $WF$.}

\Prf (a) The forward implication is clear since
$W$ is finite for $r<9$, and for $r=9$, $W$ stabilizes $-K_X$,
so assume $r>9$, or $r=9$ but $F\ne -lK_X$.
Since $p_1,\ldots,p_r$ are strongly nonspecial, if $H$ is the
class of an effective divisor, so is $wH$, for every $w\in W$.
Thus $wF\cdot H=F\cdot (w^{-1}H)\ge 0$, whenever $H$ is the class of an 
effective divisor; i.e., $wF$ is numerically effective for every $w\in W$.
Now, $F^2\ge 0$ (see, e.g., Proposition 4
of [\mtnwest]); we will first consider the case that
$F^2>0$. Then, by the index theorem, 
the subgroup $F^\perp\subset \hbox{Cl}(X)$ orthogonal to $F$
is negative definite, so the stabilizer of $F$ in $W$ is finite.
Therefore, $WF$ is infinite if $W$ is, which it is for $r\ge 9$.

Now suppose $F^2=0$. Since $e_0\cdot wF\ge 0$ for every
$w\in W$, there is a particular $w$ such that
$e_0\cdot wF$ is minimal. Let us write $wF=b_0e_0-b_1e_1-\cdots-b_re_r$ 
for some integers $b_i$. Reflections through the
roots $\rho_i$, $i>0$, just permute the coefficients
$b_1,\ldots,b_r$, so we may assume that $b_1\ge b_2\ge \cdots \ge b_r$.
In this case, if $\rho_0\cdot wF<0$, then $s_0wF\cdot e_0<wF\cdot e_0$,
contrary to assumption, so we have $wF\cdot\rho_i\ge 0$ for every
$i\ge 0$. It is not hard to show that this implies
that $wF$ is a nonnegative integer
linear combination of the classes $H_0=e_0$, $H_1=e_0-e_1$,
$H_2=2e_0-e_1-e_2$, $H_i=3e_0-e_1-\cdots-e_i$, $2<i\le r$;
i.e., $wF=\sum_ih_iH_i$ with $h_i\ge 0$. If $h_i>0$ for some $i>9$,
let $D=\sum_{i\le 9}h_iH_i+\sum_{i>9}h_iH_9$. 
Then $D^2> 0$ and $D$ is the sum of $wF$ and nonnegative
multiples of $e_{10},\ldots$, all of 
which are orthogonal to $D$, so $D$ is numerically effective. 
Thus by a previous case $W'D$ is infinite, where $W'$ is the subgroup of
$W$ generated by $s_0,\ldots,s_8$. But $W'$ stabilizes $e_{10},\ldots$,
so also $W'F$ and thus $WF$ are infinite.

So suppose that $h_i=0$ for all $i>9$. Then using 
$0\ne wF=\sum_{i\le 9}h_iH_i$
and $F^2=0$, it is easy to check that either $wF=h_1H_1$
or $wF=h_9H_9$. Since $H_9=-K_X$, if $wF=h_9H_9$, then $r>9$
by hypothesis. Let $W'$ now denote the subgroup 
generated by $s_0,\ldots,s_9$; it suffices to show $W'F$ is infinite.
I.e., it suffices to consider the case $r=10$. But if $r=10$, 
then $H_9=-K_X+e_{10}$. As is well known, $W$ fixes $K_X$
while $We_{10}$ is infinite (indeed, $We_{10}$ is the set of all
classes of $(-1)$-curves on $X$), so $WH_9$ must also be
infinite.

Finally we check that $WH_1$ is infinite.
First, $\rho=2e_0-e_4-\cdots-e_9$ is in $W\rho_1$, so reflection 
$s$ by $\rho$ is in $W$, and it is easy to check explicitly
that the composition $s_0s$ generates a cyclic subgroup $W''$
of $W$ such that $W'' H_1$ is infinite.

(b) If $wF-e_0$ is the class of an effective divisor,
then numerical effectivity of $wF$ implies
that $wF\cdot (wF-e_0)\ge 0$. Thus it suffices to show
$wF\cdot (wF-e_0)< 0$, or equivalently
$F^2< wF\cdot e_0$, for all but finitely many $wF\in WF$.
In fact, for any integer $N$ it is true that
$N<e_0\cdot wF$, for all but finitely many $wF\in WF$.
For suppose for each $D$ in an infinite subset
$V\subset WF$ we had $e_0\cdot D\le N$. Then,
writing each $D$ as $D=b_0e_0-b_1e_1-\cdots-b_re_r$ 
for integers $b_i$ depending on $D$, we would
have infinitely many integer solutions $b_0,\ldots,b_r$ to 
$F^2=b_0^2-b_1^2-\cdots-b_r^2$ with
$0\le b_0=e_0\cdot D\le N$, which is impossible.
\qed

The next result applies \ir{mycamp} to give a 
maximal rank criterion.

\prclm{Lemma}{nonunilem}{With $X$ as in \ir{Worbits}, 
let $G$ be the class of
an effective, numerically effective divisor. 
If $w\in W$ is such that there exists an $i>0$ with
$G^2<e_i\cdot w(G)$ and $G^2<(e_0-e_i)\cdot w(G)$, then
$\mu:\Gamma(w\C G)\otimes\Gamma(e_0)\to \Gamma(e_0+G)$ 
is injective, and so has maximal rank.}

\Prf By $G^2<e_i\cdot w(G)$ we have $(w(G))^2<e_i\cdot w(G)$,
but, since $w$ preserves the monoid of classes of effective divisors,
$wG$ is numerically effective, so $wG-e_i$ is not the class 
of an effective divisor; thus $q_i=h^0(X, wG-e_i)=0$. 
Similarly, $(wG-(e_0-e_i))\cdot w(G)<0$ implies
$l_i=h^0(X, wG-(e_0-e_i))=0$. Hence \ir{mycamp}
implies $\hbox{ker }\mu=0$. 
\qed

We now obtain an asymptotic result. (Given a numerically 
effective class $G$ on $X$, $Z_G$ will denote
$Z_G=m_1p_1+\cdots+m_rp_r$,
where $m_i=e_i\cdot G$.)

\prclm{Theorem}{nonunithm}{With $X$ as in \ir{Worbits},
let $G$ be the class of an effective, numerically 
effective divisor such that $h^1(X,\C G)=0$.
Then, for each $w\in W$, $I(Z_{wG})$
has the maximal rank property for all but finitely 
many elements of $\{Z_{wG}|w\in W\}$.}

\Prf Since $G$ is the class of an effective 
divisor, so is $wG$ for every $w\in W$, but, 
for all but finitely many $wG\in WG$,
$wG-e_0$ is not, by \ir{Worbits}.
Thus $\alpha(I(Z_{wG}))=wG\cdot e_0$
for all but finitely many $wG\in WG$.
On the other hand,
$h^1(X,\C G)=0$ (and hence $h^1(X,wG)=0$), so,
for all but finitely many $wG\in WG$,
the regularity of $I(Z_{wG})$ is at most $\alpha(I(Z_{wG}))+1$.
Therefore, $\mu_t(I(Z_{wG}))$ has maximal rank except possibly for
$t=\alpha(I(Z_{wG}))$; since $\mu_{\alpha(I(Z_{wG}))}$ has
maximal rank if and only if
$\mu:\Gamma(w\C G)\otimes\Gamma(e_0)\to \Gamma(e_0+wG)$
does, we turn our attention to the latter.

There are clearly only finitely many integer solutions 
$d,b_1,\cdots,b_r$ to $G^2=d^2-b_1^2-\cdots-b_r^2$
with $\{b_i: 0<i\}$ bounded. Thus the number of elements in the
orbit $WG$ with 
$\hbox{max}_{0<i}(wG\cdot e_i)\le G^2$ is finite. 
Thus it is enough by \ir{nonunilem} to show 
for each $i$ that 
$wG\cdot (e_0-e_i) > G^2$ occurs for all but finitely many 
$wG\in WG$.

We fix $i>0$; then there are 
only finitely many integer solutions 
$d,b_1,\cdots,b_r$ to $G^2=d^2-b_1^2-\cdots-b_r^2$
with $\{b_j: 0<j\ne i\}$ bounded (because then $d^2-b_i^2$
takes on only a finite set of values, which factor only a 
finite number of ways). Thus for all but finitely many 
$wG\in WG$ we can choose $0<j_w\ne i$ such that 
$wG\cdot e_{j_w}> G^2$. Now write $e_0-e_i$ as 
$(e_0-e_i-e_{j_w})+e_{j_w}$. Thus $wG\cdot (e_0-e_i)=
wG\cdot ((e_0-e_{j_w}-e_i) +e_{j_w})\ge
wG\cdot e_{j_w}> G^2$ holds 
for all but finitely many $wG\in WG$.
\qed

To apply \ir{nonunithm},
one needs examples of classes $G$ of numerically effective, effective,
and regular (i.e., $h^1=0$) divisors 
on a blowing up $X$ of \pr2 at strongly nonspecial points.
It is easy to give examples: Given such an $X$, if $m_i\ge 0$, 
then for $d$ sufficiently large (say $d>\sum_im_i$),
$G=de_0-m_1e_1-\cdots-m_re_r$ is such a class.

Alternatively, let $X$ be the blowing up of points $p_1,\ldots,p_r$
which are independent generic over the prime field.
If $-K_X\cdot G\ge 0$, then $G$ is effective,
numerically effective, and regular
if and only if $G$ is in the $W$-orbit of the nonnegative
subsemigroup $S$ of $\hbox{Cl}(X)$ generated by
$\{H_0=e_0,H_1=e_0-e_1,H_2=2e_0-e_1-e_2, H_3=3e_0-e_1-e_2-e_3,
H_4=3e_0-e_1-e_2-e_3-e_4,\ldots\}$. The proof is to
specialize $p_1,\ldots,p_r$ to a cubic, then use semicontinuity
and results of [\trans] (also see [\antican]). 

When $r\ge 9$, $W$ has a particularly tractable subgroup
for which a more explicit result analogous to 
\ir{nonunithm} can be stated (when $r<9$, $W$ is finite
and hence \ir{nonunithm} is trivial).
So assume that $p_1,\ldots,p_r\in \pr2$ are 
independent generic over the prime field
with $r\ge 9$. Let $T$ be the subgroup
of $\hbox{Cl}(X)$ generated by the roots $\rho_1,\ldots,\rho_8$.
Then, given any $v\in T$, it turns out that $v\mapsto \tau_v$
defines an injective homomorphism $T\to W$, where we define $\tau_v$
via $\tau_v(G)=G+(G\cdot H_9)v-(1/2)(2G\cdot v+(G\cdot H_9)v^2)H_9$.
If $G$ is in $S$ with $-K_X\cdot G>0$, then as above $G$ is effective,
numerically effective and regular, so, as the proof of
\ir{nonunithm} shows, $I(Z_{\tau_v(G)})$ has the maximal 
rank property for each $v\in T$ such that $G^2<e_0\cdot \tau_v(G)$,
$G^2<e_1\cdot \tau_v(G)$ and $G^2<(e_0-e_1)\cdot \tau_v(G)$.
But $T$ is negative definite and $G\cdot H_9\ge -G\cdot K_X>0$,
so substituting our expression for $\tau_v(G)$
into $e_0\cdot \tau_v(G)$, $e_1\cdot \tau_v(G)$ and $(e_0-e_1)\cdot \tau_v(G)$,
we see $G^2<e_0\cdot \tau_v(G)$, $G^2<e_1\cdot \tau_v(G)$ 
and $G^2<(e_0-e_1)\cdot \tau_v(G)$
hold for all but finitely many $v\in T$.
(In fact, we can be explicit here: these conditions 
and therefore the maximal rank property for $I(Z_{\tau_v(G)})$ hold
if $\sqrt{-v^2}> 2+\sqrt{24(G\cdot e_0)/(G\cdot H_9)}+2G^2/(G\cdot H_9)$.)

Although the foregoing paragraph provides a fairly
easy method of generating examples 
$Z=m_1p_1+\cdots+m_rp_r$ for which $I(Z)$  has the maximal rank
property, it is also nice to have an explicit criterion
in terms of the coefficients $m_i$
for the maximal rank property to hold. We give such a criterion
when $r=9$ in the next example.

\rem{Example}{explicit} Let $m_1\ge \cdots\ge m_9\ge 0$.
Here we show for general points
$p_i$ that $I(Z)$ has the maximal rank property
for $Z=m_1p_1+\cdots+m_9p_9$, if $m_1=m_9$ or if
$m_9\ge 20(m_1-m_9+1)^2$ and
$m_1+\cdots+m_9\not\equiv 2\hbox{ (mod 3)}$.

If $m_1=m_9$, which is to say that $Z$ is uniform, then it 
follows from \ir{GIGCnine} that $I(Z)$ has the maximal 
rank property, so assume that $Z$ is not uniform.
Let $X$ be the blowing up of the points $p_i$
and let $e_0,e_1,\ldots,e_9$ be the corresponding 
exceptional configuration. Since $X$ is obtained by blowing up
9 general points, $-K_X$ is numerically effective, 
so $-F_{\alpha(Z)}\cdot K_X\ge 0$, where
$F_t(Z)=te_0-m_1e_1-\cdots-m_9e_9$ (and we write just $F_t$
when the $Z$ being referred to is unambiguous). Since $F_t$
is not uniform, we have in fact that $-F_{\alpha(Z)}\cdot K_X>0$.
Thus $\alpha(Z)\ge d$, where $d$ is the largest integer
which is at most $1+(m_1+\cdots+m_9)/3$. Moreover,
$F_d^2>0$, so by Riemann--Roch $F_d$ is the class of an effective
divisor, hence actually $\alpha(Z)=d$. (To see $F_d^2>0$,
let $H=(d-3m_9)e_0-(m_1-m_9)e_1-\cdots-(m_8-m_9)e_8$,
so $F_d=H-m_9K_X$. Then $F_d^2=H^2-2m_9H\cdot K_X
\ge 2m_9-8(m_1-m_9)^2$, but by assumption we have
$m_9\ge 20(m_1-m_9+1)^2$. Note that by the same reasoning
$H-lK_X$ is the class of an effective divisor 
whenever $l\ge 4(m_1-m_9)^2$.)

Now we check that
$\alpha(Z)=\beta(Z)$; i.e., that $F_d$ is numerically
effective. As we noted above, $H-lK_X$ is effective
for $l=4(m_1-m_9)^2$. Because the points $p_i$
are general, the only curves which
could occur as fixed components of $|H-lK_X|$
of negative self-intersection are ($-1$)-curves, 
and if ${\tt E}$ is such a component, then
$E\cdot (H-lK_X)<0$. In particular, $E$ is not $e_i$ for any 
$i$, so $E\cdot e_0>0$, hence ${\tt E}$ occurs with multiplicity 
at most $e_0\cdot (H-lK_X)=d-3(m_9-l)$. Therefore,
since $-E\cdot K_X=1$, we will have
$E\cdot (H-tK_X)\ge 0$ if $t-l\ge d-3m_9+3l$,
but $d-3m_9\le 1+((m_1-m_9)+\cdots+(m_8-m_9))/3\le 4(m_1-m_9)^2$
so $m_9-l\ge 20(m_1-m_9+1)^2-4(m_1-m_9)^2\ge
16(m_1-m_9)^2\ge (d-3m_9)+3l$, so
$F_d=H-m_9K_X$ meets every such $E$ nonnegatively,
so $F_d$ is numerically effective and $h^1(X, F_d)=0$ 
by [\ir{antican}].

By definition of $d$, it is easy to check $1\le -F_d\cdot K_X\le 3$. 
Suppose $1=-F_d\cdot K_X$. Then, keeping in mind that $-K_X$ is
numerically effective, $-K_X\cdot (F_d-e_1)=0$, so
$F_d-e_1$, not being uniform but being in $K_X^\perp$, cannot be the class
of an effective divisor. Likewise,  $-K_X\cdot (F_d-(e_0-e_1))<0$, so
$F_d-(e_0-e_1)$ also cannot be the class of an effective divisor. Thus,
as in the proofs of \ir{nonunilem} and 
\ir{nonunithm}, $I(Z)$ has the maximal rank property.

Now consider the case $3=-F_d\cdot K_X$. It is easy to check
that $F_d(Z)-e_1=F_d(Z')$, where $Z'=(m_1+1)p_1+m_2p_2+\cdots+m_9p_9$,
and that  $F_d(Z)-(e_0-e_1)=F_{d-1}(Z'')$, where $Z''=
(m_1-1)p_1+m_2p_2+\cdots+m_9p_9$. 
Reasoning as above
shows that $F_d(Z')$ and $F_{d-1}(Z'')$ are classes of effective,
numerically effective divisors meeting $-K_X$ positively (and hence
$h^1(X,F_d(Z'))=0=h^1(X,F_{d-1}(Z''))$ by [\antican])
when $m_9\ge 20(m_1-m_9+1)^2$.
Arguing as above and using \ir{Firstcampcor}(c) now shows that 
$I(Z)$ has the maximal rank property.

\vskip\baselineskip

To view \ir{nonunithm} from a different perspective, 
given any distinct 
points $p_1,\ldots,p_r\in\pr2$, we define an
equivalence relation on the set of all fat point
subschemes $m_1p_1+\cdots+m_rp_r$:
we say $Z=\sum_im_ip_i$ and $Z'=\sum_im_i'p_i$ are
{\it Cremona equivalent\/} if, with respect
to the usual exceptional configuration $e_0,\ldots,e_r$
on the blow up $X$ of \pr2 at the points $p_i$,
we have $w(\beta(I(Z))e_0 -m_1e_1-\cdots-m_re_r)
=\beta(I(Z'))e_0 -m_1'e_1-\cdots-m_r'e_r$
for some $w\in W$.

When $p_1,\ldots,p_r$ are strongly nonspecial,
what it means for fat point subschemes $Z=\sum_im_ip_i$
and $Z'=\sum_im'_ip_i$ to be
Cremona equivalent  is that their associated 
linear systems in their respective degrees $\beta$
are in some sense geometrically the same 
(regarded as complete linear systems on
the blown up surface). Much else can be different,
but by the next result most (i.e., all but finitely many)
of the representatives in the equivalence class
of an expectedly good fat point subscheme 
will have ideals with the maximal rank property.

\prclm{Corollary}{applythm}{Let $p_1,\ldots,p_r\in\pr2$
be strongly nonspecial and let $Z=\sum_im_ip_i$ be
expectedly good. Then for all but finitely
many fat point subschemes $Z'=\sum_im'_ip_i$
Cremona equivalent to $Z$, $\alpha(I(Z'))=\beta(I(Z'))$ and
$I(Z')$ has the maximal rank property.}

\Prf Let $X$ and $e_0,\ldots,e_r$ be as usual and let $G=
\beta(I(Z))e_0 -m_1e_1-\cdots-m_re_r$. By definition,
every $Z'$ Cremona equivalent
to $Z$ is of the form $Z_{wG}$ for some $w\in W$. However,
not every $Z_{wG}$ need be Cremona equivalent to $Z$,
since $\beta(I(Z_{wG}))$ could be less than $wG\cdot e_0$;
in fact, from the definition we see that $Z_{wG}$ and $Z$ are 
Cremona equivalent exactly when $\beta(I(Z_{wG}))=wG\cdot e_0$.
We will check that this occurs for all but finitely many of
$\{Z_{wG}|w\in W\}$.

First, we have $\alpha(I(Z_{wG}))\le \beta(I(Z_{wG}))\le wG\cdot e_0$. 
By \ir{Worbits}(b), $wG-e_0$ is the class of an effective divisor
for at most finitely many elements $wG$ of $WG$.
This implies that $wG\cdot e_0=\alpha(I(Z_{wG}))$
and hence that $\alpha(I(Z_{wG}))=\beta(I(Z_{wG}))=wG\cdot e_0$,
for all but finitely many elements $wG$ of $WG$.
Thus all but finitely many of
$\{Z_{wG}|w\in W\}$ are Cremona equivalent to $Z$;
i.e., up to finite sets, the Cremona equivalence
class of $Z$ is $\{Z_{wG}|w\in W\}$, so \ir{applythm}
follows from \ir{nonunithm}.
\qed

We end this section applying \ir{Firstcampcor}
to uniform fat point subschemes.
In particular, in \ir{mycampcorb} and \ir{mycampcorc} we 
obtain some evidence for \ir{VIGC}, based on
the following version of \ir{Firstcampcor},
for uniform fat point ideals at 10 or 
more expectedly good points.

\prclm{Corollary}{mycampcor}{Let $p_1,\ldots,p_r$ 
be $r\ge 10$ distinct expectedly good points of \pr2, 
let $e_0,\ldots,e_r$ be the corresponding 
exceptional configuration, and let $I=I(mp_1+\cdots+mp_r)$
with $m>0$ be 
a fat points ideal. Let $F$ denote $\alpha(I)e_0-me_1-\cdots-me_r$
and define $\mu$, $l_1$, $h$ and $q_1$ as in \ir{mycamp}.
Then we have:
\item{(a)} $\alpha(I)=\beta(I)$ unless $h=1$ (in which case
$\mu$ clearly has maximal rank);
\item{(b)} the maximal rank property for $I$ holds if and only if
$\mu$ has maximal rank;
\item{(c)} $\mu$ has maximal rank if and only if
$\hbox{max}(0,2h-\alpha(I)-2) 
=\hbox{dim ker $\mu$}$; 
\item{(d)} $l_1\le \hbox{max}(0,2h-\alpha(I)-2) 
\le \hbox{dim ker $\mu$} \le l_1+q_1$; and
\item{(e)} $\hbox{max}(0,2h-\alpha(I)-2) = l_1+q_1$ 
unless $l_1=0$ and $q_1>0$.}

\Prf (a)  First we show that any uniform class 
$G=de_0-m(e_1+\cdots+e_r)$ with $m>0$ which is the class 
of an effective divisor is numerically effective
(in particular, $F$ is numerically effective).
Recall that on the blow up $X$ of expectedly good points
the only prime divisors of 
negative self-intersection are the exceptional curves
(that is, the smooth rational curves with self-intersection $-1$,
each of which thus meets $-K_X$ once). 
Now note that $d>3m$; otherwise,
$-mK_X=G+(3m-d)e_0$ is the class of an effective 
divisor with negative self-intersection meeting positively every
prime divisor of negative self-intersection, which is absurd.
But $d>3m$ means that $G$ is the class of an effective
divisor meeting every prime divisor of 
negative self-intersection positively.
Thus $G$ is numerically effective. (We also note two similarly
proved facts that we will need below: since $G-(e_0-e_1)=
(d-3m-1)e_0+e_1-mK_X$, if $G-(e_0-e_1)$ is the class of
an effective divisor, it too is numerically effective; and
$-mK_X+(d-3m-1)e_0+(e_0-e_1)=G-e_1$ so $G-e_1$
meets every exceptional curve nonnegatively, hence
$G-e_1$ is also numerically effective 
if it is effective.)

Thus $F$ is the class of an effective and numerically effective 
divisor. If $h^0(X, F)>1$, we must show that $|F|$ is free. More generally,
let ${\tt D}$ be any effective and numerically effective divisor 
on $X$ with $h^0(X,\C D)>1$. We will show that $|{\tt D}|$ is fixed
component free. Consider a Zariski decomposition $D=H+N$, where
the class of the free part of $|D|$ is $H$ and the class of
the fixed part {\tt N} is $N$. Suppose {\tt E} is an exceptional curve
which occurs as a component of {\tt N}; then $E\cdot H=0$ (else 
$h^0(X, \C D)=h^0(X, \C H)$ is impossible by Riemann--Roch), so
$0\le D\cdot E=N\cdot E$. Suppose $E\cdot C>0$ for some 
other component ${\tt C}\ne {\tt E}$ of {\tt N}.
Either {\tt C} is numerically effective or it is exceptional,
but $h^1(X, \C C)=0$ either way, so we have an exact sequence 
$0\to H^0(X, C)\to H^0(X, C+E)\to H^0({\tt E}, 
\C O_{\hbox{\eighttt E}}(C+E))\to 0$ from which the contradiction
$1=h^0(X, N)\ge h^0(X, C+E)>1$ follows. Thus {\tt E} is disjoint
from the other components of {\tt N}, and hence
$0\le D\cdot E=N\cdot E<0$. This contradiction shows 
that no exceptional curve
is a component of {\tt N}. Therefore, $N$ is numerically effective.
Thus $h^1(X,N)=0$, so $1=h^0(X, N)=1+(N^2-K_X\cdot N)/2$, 
which implies $(N^2-K_X\cdot N)/2=0$. But $1+(H^2-K_X\cdot H)/2=h^0(X, H)=
h^0(X, H+N)=1+(H^2-K_X\cdot H)/2 + H\cdot N + (N^2-K_X\cdot N)/2$, 
which implies $H\cdot N=0$. If $H^2>0$, then by 
the index theorem the subgroup
of $\hbox{Cl}(X)$ perpendicular to $H$ is negative definite; since 
$N^2\ge 0$, we must have $N=0$. Similarly, if $H^2=0$ but $N^2>0$,
then $H\cdot N=0$ implies $H=0$. Thus $H^2=0$ implies $N^2=0$
and so also $N\cdot K_X=0$. Moreover, since $1<h^0(X,H)$,
$H^2=0$ implies $-K_X\cdot H>0$. 
Now, the points $p_i$ are expectedly good, 
hence strongly nonspecial, so, as in the proof of \ir{Worbits}(a),  
$wN$ is, for some $w\in W$, a nonnegative integer
linear combination of the classes $H_i$, $0\le i\le r$.
Since $N^2=-K_X\cdot N=0$, the only possibility is that
$wN$ is a nonnegative multiple of $H_9$. If $N\ne 0$, then 
we get the contradiction: $0=H\cdot N=
wH\cdot H_9 \ge wH\cdot H_r=wH\cdot w(-K_X) =-H\cdot K_X >0$.

(b) Clearly, for $t<\alpha$ we have $I_t=0$,
so $I_t\otimes R_1\to I_{t+1}$ has maximal rank.
But the regularity of $I$ is at most $\alpha+1$
since $F$ is numerically effective and our points are
expectedly good, so $I_t\otimes R_1\to I_{t+1}$ has 
maximal rank for $t>\alpha$ by \ir{myomega}.
Thus $I$ has the maximal rank property if and only
if $\mu:I_\alpha\otimes R_1\to I_{\alpha+1}$
has maximal rank.

(c) Since $p_1,\ldots,p_r$ are expectedly good and $F$
is numerically effective, it follows that 
$h^1(X,F)=0$ so \ir{Firstcampcor}(a) implies the result.

(d) \ir{Firstcampcor}(b) gives $\hbox{max}(0,2h-\alpha(I)-2) 
\le \hbox{dim ker $\mu$}\le l_1+q_1$.

If $h^0(X,F-(e_0-e_1))=0$, 
then $l_1\le \hbox{max}(0,2h-\alpha(I)-2)$ is clear, 
so suppose $h^0(X,F-(e_0-e_1))>0$.
Thus $F-(e_0-e_1)$ is the class of an effective
divisor, hence it is numerically effective,
so $h^1(X,F-(e_0-e_1))=0$. As in the proof of 
\ir{Firstcampcor}(c) we have $h=l_1+\alpha(I)+1-m$, 
so $(2h-\alpha(I)-2)-l_1=
l_1+\alpha(I)-2m$. But $\alpha(I)^2-rm^2=F^2\ge 0$ implies
$\alpha(I)-2m>0$, and we now see $2h-\alpha(I)-2>l_1$,
which implies $l_1\le \hbox{max}(0,2h-\alpha(I)-2)$.

(e) If $q_1=0$, then (d) implies the result, so let
$q_1>0$. If also $l_1>0$, then $F-e_1$ and $F-(e_0-e_1)$
are classes of effective divisors, hence (as we saw above)
numerically effective, so $h^1(X,F-e_1)=0=h^1(X, F-(e_0-e_1))$,
so the result follows by \ir{Firstcampcor}(c).
\qed

\rem{Remark}{campcomp} Whereas the bound $\hbox{max}(0,2h-\alpha(I)-2) 
\le \hbox{dim ker $\mu$}$ in \ir{mycampcor}
is in fact exactly what one obtains
from [Cam], the upper bound $\hbox{dim ker $\mu$}\le l_1+q_1$
is always at least as good as Campanella's
(which is always either $h-1$ or $h-2$), and except in extremal cases
(i.e., $h\le 2$ or $\alpha\le h\le \alpha+1$) it is better.

\vskip\baselineskip

Assuming expectedly good points, computer runs suggest 
that $\hbox{max}(0,2h-\alpha(I)-2)$ equals $l_1+q_1$
fairly often, possibly for infinitely many $m$ for each
$r>9$ which is not an even square. The next two
corollaries verify this possibility
for some special values of $r$.

\prclm{Corollary}{mycampcorb}{Using the notation and
hypotheses of \ir{mycampcor}, $\mu$ has maximal rank
for infinitely many $m$ whenever $r+i$ is an odd square
for some $i\in\{-3, -2, -1, 0,1,2,3,4\}$.}

\Prf First assume $r+i$ is an odd square for some $i\in\{0,1,2,3,4\}$;
then it is not hard to see that there is an odd integer 
$2t+1\in [\sqrt{r},\sqrt{r}+2/\sqrt{r})$. By \ir{mycampcor},
$\mu$ has maximal rank whenever $q_1=0$, so for the given $r$ 
it suffices to check that $q_1=0$ for infinitely many $m$. 

By the proof of \ir{mycampcor}, 
$F-e_1$ is numerically effective 
whenever it is effective.
Taking cohomology of $0\to \C O_X(F-e_1)
\to \C F\to \C F\otimes \C O_{{\hbox{\eighttt E}}_1}\to 0$ we see
the restriction map $H^0(X,\C F)\to H^0(E_1,\C F\otimes
\C O_{{\hbox{\eighttt E}}_1})$
always has maximal rank. Thus $h^0(X,\C F)\le 
h^0({\tt E}_1,\C F\otimes\C O_{{\hbox{\eighttt E}}_1})$
implies that $q_1=0$.

To apply this, note that $h^0({\tt E}_1,\C F\otimes
\C O_{{\hbox{\eighttt E}}_1})=m+1$
and, since $F$ is numerically effective 
by the proof of \ir{mycampcor}, that
$h^0(X,\C F)=\binom{\alpha(I)+2}{2}-r\binom{m+1}{2}$. 
From this we obtain the following criterion: $q_1=0$ for each 
$m$ for which $0<\binom{x+2}{2}-r\binom{m+1}{2}\le m+1$
has a positive integer solution $x$.

Now, we know $b^2-rm^2=\epsilon$ has infinitely many positive integer
solutions $(b,m)$, where we take $\epsilon=0$ if $r$ is a square
and we take $\epsilon=1$ otherwise (in 
which case we have Pell's equation).
Substituting $b+t-1$ in our criterion for $x$ and simplifying gives
$-t^2-t-\epsilon<(2t+1)b-rm\le 2m+2-t^2-t-\epsilon$.

Since $b\ge \sqrt{r}m$ and $2t+1\ge \sqrt{r}$, we clearly have 
$-t^2-t-\epsilon<(2t+1)b-rm$. 
Since $b<\sqrt{r}m+1$, we see $(2t+1)b-rm<(2t+1)(\sqrt{r}m+1)-rm$;
i.e., $(2t+1)b-rm$ is bounded above by a 
linear function of $m$. Using $2t+1<\sqrt{r}+2/\sqrt{r}$
shows that the coefficient of $m$ in this linear function
is less than 2, so for $m$ sufficiently
large we have $(2t+1)b-rm<2m+2-t^2-t-\epsilon$.

Now assume $r+i$ is an odd square for some $i\in\{-1,-2,-3\}$;
then it is not hard to see that there is an odd integer 
$2t-1\in (\sqrt{r}-2/\sqrt{r},\sqrt{r})$. By \ir{mycampcor}
it suffices to check that $l_1>0$ for infinitely many $m$,
so this time we use the fact that
$F-(e_0-e_1)$ is numerically effective 
whenever it is effective. From the proof of \ir{mycampcor}(d), 
we have $l_1=\hbox{max}(0,h-\alpha+m-1)$.
Thus $h-\alpha+m-1>0$ implies $l_1>0$, which gives us
the following criterion: $l_1>0$ for each 
$m$ for which $x <\binom{x+1}{2}-r\binom{m+1}{2} +m\le x+ m$
has a positive integer solution $x$. (If $x$ is a solution, then
$x=\alpha+1$. In particular, the second inequality fails for
$x > \alpha+1$, while the first fails for
$x<\alpha$.) Simplifying gives
$x +(r-2)m< x^2-rm^2\le x+rm$, and as above, $b^2-rm^2=1$ 
has infinitely many positive integer solutions $(b,m)$.
Substituting $b+t$ in for $x$ and simplifying gives
$2m-(t-t^2-1)>rm-(2t-1)b\ge -t+t^2+1$.

Since $b\ge \sqrt{r}m$, we have $rm-(2t-1)b\le rm-(2t-1)\sqrt{r}m$,
so $rm-(2t-1)b$ is bounded above by a linear function of $m$,
and using $2t-1>\sqrt{r}-2/\sqrt{r}$ shows 
that the coefficient of $m$ in this linear function is less than two.
It now follows that our criterion's first inequality holds for all
sufficiently large $m$. 
For the other inequality, using $b<\sqrt{r}m+1$ shows
$rm-(2t-1)b$ is strictly bounded below by a linear function of $m$,
and now using $2t-1<\sqrt{r}$ shows the coefficient
of $m$ in this linear function is positive.
Thus the second inequality also holds for all
sufficiently large $m$. \qed

\ir{mycampcorb} gives a partial answer
to the question of for which $r$ do our bounds
infinitely often force $\mu$ to have maximal rank.
An interesting side remark here is that in fact the bounds force
$\mu$ to have maximal rank for all but finitely many $m$ when $r$ 
is an odd square, whereas for $r$ an even square our bounds never
force maximal rank. 
A slightly different approach (and further easy
variations of it) gives additional examples, such as $r=13$.

\prclm{Corollary}{mycampcorc}{Using the notation and
hypotheses of \ir{mycampcor}, $\mu$ has maximal rank
for infinitely many $m$ whenever $r=(ca)^2+4c^2>9$ 
for positive odd integers $a$ and $c$.}

\Prf Here we use the criterion developed in the proof of
\ir{mycampcorb} involving $q_1=0$.
So substitute $x=\sqrt{r}m+t-1$ into
$0<\binom{x+2}{2}-r\binom{m+1}{2}\le m+1$. For $t$ real, this has
solutions for all sufficiently large integers $m>0$ if 
$t$ is in the interval 
$[(\sqrt{r}-1)/2, (\sqrt{r}-1)/2 +1/\sqrt{r})$,
and so for $m$ sufficiently large
$0<\binom{x+2}{2}-r\binom{m+1}{2}<m+1$ has a positive integer
solution $x$ if the interval
$[\sqrt{r}m+(\sqrt{r}-1)/2, \sqrt{r}m+(\sqrt{r}-1)/2 +1/\sqrt{r})$
contains an integer. After simplifying, this is equivalent 
to finding an integer $\lambda$ such that
$$0\le {{2\lambda+1}\over{2m+1}}-\sqrt{r}< {{2}\over{\sqrt{r}(2m+1)}}.$$
It is well known that $\sqrt{r}$ 
has infinitely many rational approximations
$p/q$ accurate to order $1/q^2$ if $r$ is 
not a square. The problem here is to ensure
in addition that $p$ and $q$ are odd with $p/q>\sqrt{r}$. Whether
this also is known we do not know, but it can at least be verified
in certain cases. For example, consider the continued fraction expansion
$$\sqrt{(ca)^2+4c^2} = ca + {2c\over\displaystyle a+
                     {\strut 1\over\displaystyle a+\cdots}}.$$
Taking successive convergents (see [Bk] for background
on continued fractions)
gives a sequence $\{c_i\}$ of rational approximations
which for $i\equiv 2 (\hbox{mod }6)$ is a ratio $p/q>\sqrt{r}$
of odd integers $p$ and $q$. Moreover, the general theory
of continued fractions implies that each convergent $p/q$ is accurate
to order $1/q^2$. Thus for $r=(ca)^2+4c^2>9$ expectedly good points, 
$\mu$ has maximal rank for infinitely many $m$.
\qed

\irrnSection{Results for $r<10$}{smallr}
Finally, we prove complete results for arbitrary symbolic powers
of ideals of $r\le 9$ general 
points of \pr2; i.e., for ideals $I(Z)$ of 
fat point subschemes $Z=mp_1+\cdots+mp_r$
for $r\le 9$ general points $p_i$. Along the way we prove \ir{data}
and we conclude by proving \ir{mxrkfails}.

We divide our analysis into three cases, $r\le 5$, $6\le r\le 8$,
and $r=9$, with the second case requiring most of the effort but also 
being the most interesting.

\irSubsection{Five or Fewer General Points}{conic}
Let $X$ be the blow up of \pr2 at $r\le 5$ general points.
By \ir{cokfact} and \ir{recall}, we can compute
$\s(\C F,e_0)$ for an arbitrary class $F$ if
we can do so whenever $F$ is a numerically effective class.
But any five or fewer general points in the plane lie on a smooth conic,
so the results of [Cat] apply. Translating the 
results of [Cat] to the language used here and
examining what [Cat] proves, we find that 
$\C S(\C F,e_0)=0$ for any numerically 
effective class $F$. 
(In fact, [Cat] iteratively finds generators for and a resolution of
$I(Z)$ for any fat point subscheme $Z=m_1p_1+\cdots+m_tp_t$, where
$p_1,\ldots,p_t$ are distinct points of a smooth plane conic,
which includes the case of a 
uniform $Z$ supported at 5 or fewer general points
of \pr2. From our perspective, 
the key fact in [Cat], not explicitly stated there, is that
$\C S(\C F,e_0)=0$ for any numerically effective class $F$ on the blow up 
$X$ of points on a smooth conic. See [\fatpts] for 
an explicit proof and a generalization.)

Applying the foregoing to $Z=m(p_1+\cdots+p_r)$ for
$r\le 5$ general points 
$p_1,\ldots,p_r\in\pr2$ and $m>0$, we have the following.
Since $F=\beta(I(Z))e_0-m(e_1+\cdots+e_r)$ is numerically effective
and hence $\s(\C F,e_0)=0$,
we see $\mu_\beta$ is surjective and so has maximal rank;
thus the RUMRP holds for $r\le 5$. As for the UMRP,
for $r=1$ it is easy to see that 
$I(Z)_t=0$ for $t<m$ and that $te_0-me_1$ is
numerically effective for $t\ge m$. 
The former means that $\mu_t$ is 
injective for $t<m$, and by the preceding paragraph
and numerical effectivity of $te_0-me_1$ for $t\ge m$, we have 
$\C S(te_0-me_1,e_0)=0$ for $t\ge m$, 
and hence $\mu_t$ is surjective for $t\ge m$. 
Thus the UMRP holds on \pr2 for $r=1$.
For $r=4$, $I(Z)_t=0$ for $t<2m$, since $2e_0-(e_1+\cdots+e_4)$ is 
numerically effective but 
$[2e_0-(e_1+\cdots+e_4)]\cdot[te_0-m(e_1+\cdots+e_4)]<0$.
Also, $\C S(\C F_t,e_0)=0$ for 
$F_t=te_0-m(e_1+\cdots+e_4)$ with $t\ge 2m$,
since $F_t=m(2e_0-(e_1+\cdots+e_4))+(t-2m)e_0$ 
is numerically effective. 
Thus the UMRP holds on \pr2 also for $r=4$.  

To see that the UMRP on \pr2 fails for $r=2,3,5$, 
it is enough by \ir{abfail}
to find in each case a uniform abnormal class. 
But these have already been exhibited in
\ir{recall}: for $r=2$, we have
$e_0-(e_1+e_2)$; for $r=3$, there is $3e_0-2(e_1+e_2+e_3)$;
and for $r=5$, $2e_0-(e_1+\cdots+e_5)$. 
One can check, in fact, that for $r=2,3,5$,
$I(m(p_1+\cdots+p_r))$ fails to 
have the maximal rank property if and only if:
$r=2$ and $m\ge 2$; or $r=3$ or $r=5$ and $m\ge 3$.

\irSubsection{Six to Eight General Points}{sixtoeight}
\ir{subst} determines $\s(\C F,e_0)$ for any numerically effective
uniform class on a blow up $X$ of \pr2 at $6\le r\le 8$ general points.
That the RUMRP holds for $r=6$ but not for $r=7$ or 8
follows directly from \ir{subst}. That the UMRP fails 
for $6\le r\le 8$ follows from \ir{abfail}, since,
as shown in \ir{recall}, in each case
the blow up of $r$ general points supports a uniform abnormal class:
for $r=6$, $E=12e_0-5(e_1+\cdots+e_6)$ is such;
for $r=7$, $E=21e_0-8(e_1+\cdots+e_7)$ is such; and 
for $r=8$, $E=48e_0-17(e_1+\cdots+e_8)$ is such. 

\prclm{Theorem}{subst}{Let $F=F(d,m,r)$ 
be a uniform numerically effective class
on the blowing up $X$ of $6\le r\le 8$ general 
points of \pr2 (where $F(d,m,r)$
denotes $de_0-m(e_1+\cdots+e_r)$). 
\itemitem{(a)} If $r=6$, then $\r(e_0,\C F)\s(e_0,\C F)=0$.
\itemitem{(b)} If $r=7$, then $\r(e_0,\C F)\s(e_0,\C F)=0$ unless 
$F=lF(8,3,7)$ for $l\ge 3$, in which case $\s(e_0,\C F)=7$.
\itemitem{(c)} If $r=8$, then $\r(e_0,\C F)\s(e_0,\C F)=0$, unless
$F=lF(17,6,8)$ for $l\ge 9$, in which case $\s(e_0,\C F)=48$, or
unless $F=lF(17,6,8)+F(3,1,8)$ for 
$l\ge 6$, in which case $\s(e_0,\C F)=16$.}

\Prf Note that we can compute $h^0$ 
for any class $F$ on $X$, as discussed in
\ir{recall} or more generally using [\trans], keeping in mind
that any 8 or fewer general points are expectedly good [\mtnwest].

So let $6\le r\le 8$ and let $F=de_0-m(e_1+\cdots+e_r)$ be a uniform
class. If $F$ is numerically effective, then $h^1(X,F+te_0)=0$ 
for all $t\ge 0$ (by \ir{recall}), so $\C S(F+te_0,e_0)=0$
for all $t>0$ by \ir{myomega}. Thus 
we only need to consider $\delta e_0-m(e_1+\cdots+e_r)$,
where $\delta$ is the least $d$ such that 
$de_0-m(e_1+\cdots+e_r)$ is numerically effective.
We will denote this class as $F_m$, or 
just by $F$ if our meaning is clear.
Using \ir{recall} it follows that $\delta$
is the least positive integer $d$ such that: $d\ge 5m/2$ if $r=6$;
$d\ge 8m/3$ if $r=7$; or $d\ge 17m/6$ if $r=8$. 

First say $r=6$ and ${\tt E}$ is the effective divisor whose class is
$E$, where here we take $E=12e_0-5(e_1+\cdots+e_6)$. Note that
${\tt E}$ is a disjoint union of six $(-1)$-curves.
Also, if $m$ is odd, then $F_m=-K_X+(m-1)(5e_0-2(e_1+\cdots+e_6))/2$,
while $F_m=m(5e_0-2(e_1+\cdots+e_6))/2$ if $m$ is even. 
In any case, $h^2(X,F_m-e_0)=0$ by duality.

If $m$ is odd, one checks (by induction on $m$ in 
$0\to\C O_X(F_m-e_0)\to
\C O_X(F_{m+2}-e_0)\to\C O_X(\C F_{m+2}-e_0)\otimes
\C O_{{\hbox{\eighttt C}}})\to0$,
where the class of the smooth rational curve ${\tt C}$ is $F_2$)
that, suppressing the subscript, $h^1(X,F-e_0)=0$,
and hence (by \ir{myomega}) that $\C S(\C F,e_0)=0$
so suppose $m=2s$, with $s\ge 1$. For $s=2$,
$e_0\cdot(F-E+e_0)=-1$, so $h^0(X,F-E+e_0)=0$
so $\C S(F-E,e_0)=0$. For $s>2$, $F-E$
is numerically effective with odd uniform multiplicity,
so $\C S(F-E,e_0)=0$ by the preceding case.
Since $\C F\otimes\C O_{{\hbox{\eighttt E}}}=
\C O_{\hbox{\eighttt E}}$, it is easy to check that 
$\C S(\C O_{\hbox{\eighttt E}},e_0)=0$, using \ir{Mumford}(b) 
applied to the components of ${\tt E}$. 
If we now check that $h^1(X,F-E+e_0)=0$ and
$h^1(X,F-E)=0$, then we can apply \ir{Mumford}(a)
to $(0\to\Gamma(F-E)\to\Gamma(\C F)\to
\Gamma(\C F\otimes\C O_{\hbox{\eighttt E}})\to0)\otimes\Gamma(e_0)$
to obtain $\C S(\C F,e_0)=0$.
But for $s>2$, we have $h^1(X,F-E+e_0)=0$ and
$h^1(X,F-E)=0$ by \ir{recall}. 
For $s=2$, we have $F-E=K_X+e_0$
and $F-E+e_0=K_X+2e_0$; now using duality 
and descending to \pr2 we see
$h^1(X,F-E+ae_0)=h^1(\pr2,\C O_{\pr2}(-a-1))=0$ for any $a$.

We are left with the case $s=1$, thus $m=2$, but 
here $(F-(e_0-e_1))\cdot F<0$ and $(F-e_1)\cdot F<0$ 
so $l_1=q_1=0$ and, by \ir{Firstcampcor}, $\r(\C F,e_0)=0$.

Now say $r=7$ and ${\tt E}$ is the effective divisor whose
class is $E=21e_0-8(e_1+\cdots+e_7)$. This time ${\tt E}$ 
is a union of seven disjoint 
$(-1)$-curves and $F_m=s\C F_3-tK_X$, where 
$F_3=8e_0-3(e_1+\cdots+e_7)$
and the integers $s$ and $t$ are defined by taking
$m=3s+t$ such that $0\le t<3$.

For $t=2$ and any $s\ge 0$ we have 
$h^1(X,F_m-e_0)=0$ (as in the case $r=6$), 
which gives $\C S(\C F_m,e_0)=0$ by \ir{myomega}. 
For $m<9$, we have $\r(\C F_m,e_0)\s(\C F_m,e_0)=0$ 
(with, in fact, $\s(\C F_m,e_0)=0$ when $m$ is not 3 or 6) by 
computing cohomology and applying
\ir{Firstcampcor}. Similarly, for $m=10$ we have $\s(\C F_m,e_0)=0$, so
applying \ir{Mumford}(a) with $m=10$ to 
$(0\to\Gamma(F_m-E)\to\Gamma(\C F_m)\to
\Gamma(\C F_m\otimes\C O_{\hbox{\eighttt E}})\to0)\otimes\Gamma(e_0)$
shows that $\s(\C F_m\otimes
\C O_{\hbox{\eighttt E}},e_0)=0$. But $F_{3s+1}\otimes
\C O_{\hbox{\eighttt E}}=
F_{10}\otimes\C O_{\hbox{\eighttt E}}$ for any $s$,
hence $\s(\C F_{3s+1}\otimes\C O_{\hbox{\eighttt E}},e_0)=0$ for any $s$.
Checking $\s(\C F_{3s}\otimes
\C O_{\hbox{\eighttt E}},e_0)=7$ is even easier, using 
$\C F_{3s}\otimes\C O_{\hbox{\eighttt E}}=\C O_{\hbox{\eighttt E}}$.
We can now handle the remaining cases, $0\le t\le 1$
with $s\ge 3$; for these we consider
$(0\to\Gamma(F_m-E)\to
\Gamma(\C F_m)\to\Gamma(\C F_m\otimes\C O_{\hbox{\eighttt E}})\to0)\otimes
\Gamma(e_0)$,
using $m=3s+t$, $F_m-E=F_{3(s-3)+t+1}$ and 
$\C F_m\otimes\C O_{\hbox{\eighttt E}}=\C O_X(-tK_X)\otimes
\C O_{\hbox{\eighttt E}}$.
By induction, $\s(\C F_{3(s-3)+t+1},e_0)=0$, 
so by the exact sequence of
\ir{Mumford}(a) we obtain $\s(\C F_{3s+1},e_0)=0$ 
and $\s(\C F_{3s},e_0)=7$.

In conclusion, for $r=7$, $\r(\C F,e_0)\s(\C F,e_0)=0$ for all
numerically effective uniform classes $F$ except $F=lF_3$
for $l\ge3$, in which case $\s(\C F,e_0)=7$.

We now proceed to the last case, for 
which $X\to\pr2$ is a blow up of
$r=8$ general points of \pr2. Here we 
let ${\tt E}$ be the effective divisor whose class is 
$E=48e_0-17(e_1+\cdots+e_8)$; ${\tt E}$ is 
a union of eight disjoint $(-1)$-curves, each of
which under $X\to\pr2$ maps to a plane 
sextic with seven double points and a triple point.
Here we have
$F_m=sF_6-tK_X$, where 
$F_6=17e_0-6(e_1+\cdots+e_8)$ and $m=6s+t$ with
$0\le t<6$.
It follows for $s\ge 3$ that 
we have $F_{6s+t}-E=F_{6(s-3)+t+1}$.

We first need to compute $\s(\C F_m\otimes\C O_{\hbox{\eighttt E}},e_0)$
and $\r(\C F_m\otimes\C O_{\hbox{\eighttt E}},e_0)$ for $0\le t<5$.
Let $C$ be the class of any component 
${\tt C}$ among the eight components of ${\tt E}$. Then
$\s(\C F_m\otimes\C O_{\hbox{\eighttt E}},e_0)=
8\s(\C F_m\otimes\C O_{\hbox{\eighttt C}},e_0)$, 
so we restrict our attention
to ${\tt C}$. Note that $\C F_m\otimes\C O_{\hbox{\eighttt C}}
=\C O_X(-tK_X)\otimes\C O_{\hbox{\eighttt C}}=\C O_{\hbox{\eighttt C}}(t)$.

For $t=0$, clearly $\r(\C O_{\hbox{\eighttt C}},e_0)=0$
(a linear form times a nonzero
constant cannot vanish on a sextic) whence 
$\s(\C O_{\hbox{\eighttt C}},e_0)=4$, so consider
$t=1$. Then we have $\r(\C O_{\hbox{\eighttt C}}(1),e_0)=0$ 
and so $\s(\C O_{\hbox{\eighttt C}}(1),e_0)=2$:
letting $x$ and $y$ be a basis for $\Gamma(\C O_{\hbox{\eighttt C}}(1))$, 
a nontrivial element of 
$\Gamma(\C O_{\hbox{\eighttt C}}(1))\otimes\Gamma(e_0)$ which maps to 0 in 
$\Gamma(\C O_{\hbox{\eighttt C}}(1)\otimes e_0)=
\Gamma(\C O_{\hbox{\eighttt C}}(7))$ gives an equation
$xf=yg$, where $f$ and $g$ are 
restrictions to $C$ of distinct lines in \pr2.
But $f$ and $g$ have degree 6, so $xf=yg$ 
implies $f$ and $g$ have 5 zeros on ${\tt C}$ in common.
Since the image of ${\tt C}$ in \pr2 
has at most a triple point, two distinct
lines can have at most 3 points of ${\tt C}$ in common, 
contradicting there being a nontrivial
element of the kernel.

For $t=2$, both $\C R(\C O_{\hbox{\eighttt C}}(2),e_0)$ and 
$\C S(\C O_{\hbox{\eighttt C}}(2),e_0)$ vanish:
let $x$ and $y$ be as before and let 
$f,g,h$ be a basis for the restriction
of $\Gamma(e_0)$ to ${\tt C}$ such that $f$ and 
$g$ correspond to lines in \pr2
which meet at the triple point of the image of ${\tt C}$ in \pr2.
If $\r(\C O_{\hbox{\eighttt C}}(2),e_0)\ne0$, then we have an equation
$q_1f+q_2g+q_3h=0$, where $q_1,q_2,q_3$ (not all 0) lie in the span of
$\{x^2,xy,y^2\}$. Since $f$ and $g$ have exactly 3 zeros in common,
we cannot have $q_3=0$, and so $h$ 
also has a zero in common with $f$ and $g$,
which gives the contradiction that the 
restriction of $\Gamma(e_0)$ to ${\tt C}$
has a base point. Thus $\r(\C O_{\hbox{\eighttt C}}(2),e_0)=0$ 
from which we easily compute
$\s(\C O_{\hbox{\eighttt C}}(2),e_0)=0$.

For $t=3,4$ or 5, we have $\s(\C O_{\hbox{\eighttt C}}(t),e_0)=0$: 
say $t=3$ ($t=4$ or 5 are similar).
Let $x$ and $y$ be as above; thus cubics in $x$ 
and $y$ span $\Gamma(\C O_{\hbox{\eighttt C}}(3))$.
But $\Gamma(\C O_{\hbox{\eighttt C}}(1))\otimes
\Gamma(\C O_{\hbox{\eighttt C}}(2))$ surjects onto
$\Gamma(\C O_{\hbox{\eighttt C}}(3))$, and, by the previous case,
$\Gamma(\C O_{\hbox{\eighttt C}}(2))\otimes\Gamma(e_0)$ surjects onto
$\Gamma(\C O_{\hbox{\eighttt C}}(8))$, so 
$\Gamma(\C O_{\hbox{\eighttt C}}(3))\otimes\Gamma(e_0)$
and $(\Gamma(\C O_{\hbox{\eighttt C}}(1))\otimes
\Gamma(\C O_{\hbox{\eighttt C}}(2)))\otimes\Gamma(e_0)$
and $\Gamma(\C O_{\hbox{\eighttt C}}(1))\otimes
\Gamma(\C O_{\hbox{\eighttt C}}(8))$
all have the same image in 
$\Gamma(\C O_{\hbox{\eighttt C}}(9))$. Since ${\tt C}$ is rational,
we know $\Gamma(\C O_{\hbox{\eighttt C}}(1))\otimes
\Gamma(\C O_{\hbox{\eighttt C}}(8))$
surjects onto $\Gamma(\C O_{\hbox{\eighttt C}}(9))$, 
whence $\s(\C O_{\hbox{\eighttt C}}(t),e_0)=0$.

Now we are ready to consider $\s(\C F_m,e_0)$. First,
$h^1(X,F_{6s+5}-e_0)=0$ for $s\ge 0$, hence
$\C S(\C F_{6s+5},e_0)=0$ for all $s\ge 0$ by \ir{myomega}.
For $m=2, 8, \hbox{ or } 14$ we apply
\ir{Mumford} to
$(0\to\Gamma(F_m-C)\to\Gamma(\C F_m)\to
\Gamma(\C F_m\otimes\C O_{\hbox{\eighttt C}})\to0)\otimes\Gamma(e_0)$,
where we take $C$ to be, respectively,
$6e_0-3e_1-2e_2-\cdots-2e_8$,
$24e_0-9(e_1+\cdots+e_4)-8(e_5+\cdots+e_8)$ and
$42e_0-15(e_1+\cdots+e_7)-14e_8$, from which it follows that
$\s(\C F_m,e_0)$ is, respectively, 0, 1 and 0, from which we derive
that $\r(\C F_m,e_0)\s(\C F_m,e_0)=0$. 
We also check $\r(\C F_m,e_0)\s(\C F_m,e_0)=0$ for $21\le m\le 22$ and for
the remaining values
of $0\le m\le 18$ by applying \ir{Firstcampcor} (we suppress
the explicit computations).

It turns out, in fact, that 
$\C S(\C F_{6s+t},e_0)=0$ for $3\le t\le 4$ and
$0\le s\le 3$. Thus, using $\C S(\C F_{6s+5},e_0)=0$
and $\C S(\C F_{6s+t}\otimes\C O_E,e_0)=0$ for 
$3\le t\le 4$ from above, with 
$(0\to\Gamma(F_m-E)\to\Gamma(\C F_m)\to
\Gamma(\C F_m\otimes\C O_{\hbox{\eighttt E}})\to0)\otimes\Gamma(e_0)$ 
and \ir{Mumford} and induction, we conclude 
$\C S(\C F_{6s+t},e_0)=0$ for $3\le t\le 4$ and
all $s\ge 0$. In particular, we now see that either 
$\C R(\C F_{m}\otimes\C O_{\hbox{\eighttt E}},e_0)=0$
or $\C S(F_m-E,e_0)=0$ for every $m\ge 18$, and hence
that $0\to \C S(\C F_{6(s-3)+t+1},e_0)\to
\C S(\C F_{6s+t},e_0)\to\C S(\C F_{6s+t}\otimes\C E,e_0)\to0$
is exact for all $0\le t< 5$ and any $s\ge 3$.

For $s\ge 3$, we thus obtain three recursion formulas:
$\s(\C F_{6s+2},e_0)=\s(\C F_{6(s-3)+3},e_0)+
\s(\C F_{6s+2}\otimes\C O_{\hbox{\eighttt E}},e_0)$
(or $\s(\C F_{6s+2},e_0)=0$, since $\s(\C F_{6(s-3)+3},e_0)$
and $\s(\C F_{6s+2}\otimes\C O_{\hbox{\eighttt E}},e_0)$ vanish);
$\s(\C F_{6s+1},e_0)=\s(\C F_{6(s-3)+2},e_0)+16$ 
(since $\s(\C F_{6s+1}\otimes\C O_{\hbox{\eighttt E}},e_0)=16$); and
$\s(\C F_{6s},e_0)=\s(\C F_{6(s-3)+1},e_0)+32$ 
(since $\s(\C F_{6s}\otimes\C O_{\hbox{\eighttt E}},e_0)=32$).
Thus we see that $\s(\C F_{6s+2},e_0)=0$ for $s>1$,
that $\s(\C F_{6s+1},e_0)=16$ for $s>4$, and 
that $\s(\C F_{6s},e_0)=48$ for $s>7$.
It is also now easy to check 
that $\r(\C F_{6s+1},e_0)=0$ if and 
only if $s<6$, and that $\r(\C F_{6s},e_0)=0$ 
if and only if $s<9$. Thus with our results above we have that
$\r(\C F_m,e_0)\s(\C F_m,e_0)=0$ for every $m\ge 0$ except
when $m=6s$ and $s\ge 9$, in which case $\s(\C F_{m},e_0)=48$,
or when $m=6s+1$ and $s\ge 6$, in which case $\s(\C F_{m},e_0)=16$.\qed

\rem{Example}{example} Here we use our 
results to explicitly compute a resolution's modules.
Consider eight general points, each
taken with multiplicity $m=205$; thus $Z=205(p_1+\cdots+p_8)$.
Then $I(Z)$ has: $\nu_{579}=10$ generators in degree 579 (since
579 is the first degree $d$ such that $I(Z)_d\ne 0$, and we have
$\hbox{dim}_kI(Z)_{579}=10$); $\nu_{580}=201$
(since $579e_0-205(e_1+\cdots+e_8)$ has 
free part $H=51e_0-18(e_1+\cdots+e_8)$
and fixed part $N=528e_0-187(e_1+\cdots+e_8)$, and here
$\s(e_0,H+N)=\s(e_0,\C H) + (h^0(X,e_0+H+N)-h^0(X,H+e_0))
=33+168$); $\nu_{581}=208$ (since, for $d=580$, 
$H=340e_0-120(e_1+\cdots+e_8)$
and $N=240e_0-85(e_1+\cdots+e_8)$, and
$\s(e_0,H+N)=48+160$); and $\nu_{582}=16$
(since, for $d=581$, $H=581e_0-205(e_1+\cdots+e_8)$
and $N=0$, and $\s(e_0,H+N)=16+0$). Moreover, the regularity
of $I(Z)$ is 582, so there are no other generators.
(These numbers can be compared with Campanella's bounds
[Cam]: $10\le \nu_{579}\le 10$, $201\le \nu_{580}\le 210$, 
$70\le \nu_{581}\le 280$, and $0\le \nu_{582}\le 79$.)

Now, a minimal free resolution of $I(Z)$ has the form
$0\to F_1\to F_0\to I(Z)\to 0$, and
our data show that
$F_0=R^{10}[-579]\oplus R^{201}[-580]\oplus 
R^{208}[-581]\oplus R^{16}[-582]$.
We therefore now know the Hilbert function of $F_0$,
which with the Hilbert function of $I(Z)$
determines the Hilbert function of $F_1$, from which
we recover $F_1$ itself as $R^{138}[-581]\oplus 
R^{216}[-582]\oplus R^{80}[-583]$.

\irSubsection{Nine General Points}{nine}
Nine general points of \pr2 always lie on a smooth
cubic curve, so let
$X$ be obtained by blowing up nine 
distinct points $p_1,\ldots,p_9$ of a smooth
cubic curve in \pr2. As indicated in \ir{recall},
$-K_X$ is the class of a smooth elliptic curve ${\tt C}$, and
any uniform class $F=de_0-m(e_1+\cdots+e_9)$ on $X$ with $m>0$
can be written as $F=te_0-mK_X$ (with $t=d-3m$).
As in \ir{recall}, if the restriction of $\C O_X(-tK_X)$ to $C$ has infinite
order in $\hbox{Pic}({\tt C})$, define $a$ to be 0.
If the order $l$ is finite, define $a$ by 
requiring $m=al+b$ with $0\le b<l$.

\prclm{Theorem}{unifnine}{Let $F=te_0-mK_X$ with $m>0$,
where $X$ is the blowing up of $9$ distinct points of a
smooth cubic in \pr2, with $a$ defined as above.
\item{(a)} If $t<-1$, then $\C S(\C F,e_0)=0$.
\item{(b)} If $t=-1$, then $\s(\C F,e_0)=a+1$ and $\r(\C F,e_0)=0$.
\item{(c)} If $t=0$, then $\s(\C F,e_0)=3m-3a$ and $\r(\C F,e_0)=0$.
\item{(d)} If $t>0$, then $\s(\C F,e_0)=0$.
}

\Prf Part (a) follows since then $h^0(X,F+e_0)=0$,
part (b) follows since then $\s(\C F,e_0)=h^0(X,F+e_0)=a+1$,
while parts (c) and (d) follow from 
Theorem 3.2.1.2 of [\fatpts], and, in the case of (c), 
from \ir{cokfact} and \ir{recall}.\qed

\prclm{Corollary}{GIGCnine}{The UMRP (and hence
the RUMRP) holds on \pr2 for $r=9$.}

\Prf By \ir{unifnine}, we see 
$\r(\C F,e_0)\s(\C F,e_0)=0$ for any uniform class
on the blowing up $X$ of any 9 distinct points of a smooth cubic. 
Thus $I(Z)$ has the maximal 
rank property for any uniform
fat point subscheme $Z$ supported at nine general 
points of \pr2.\qed

\irSubsection{Proof of \ir{mxrkfails}}{corproof}

We close with the proof of \ir{mxrkfails}.  Let 
$p_1, \ldots,p_r$ be $r\le 9$ general points of \pr2.
Then \ir{subst} shows that $I(m(p_1+\cdots+p_r))$
fails to have the maximal rank property if:
$r=7$, $m=3l$ and $3\le l\le 7$;
or $r=8$, $m=6l$ and $9\le l\le 16$; or
$r=8$, $m=6l+1$ and $6\le l\le 13$.
Using \ir{recall} we check in each of these cases that
$\alpha(I)=\beta(I)$.

Conversely,  since $r\le 9$ and $\alpha(I)=\beta(I)$,  
$I$ has the maximal rank property if and only if 
$\mu_{\beta(I)}$ has maximal rank.
By \ir{conic}, \ir{subst} and \ir{GIGCnine},
if $\mu_{\beta(I)}$ fails to have maximal rank, 
then either $r=7$, $m=3l$ and $3\le l$,
or $r=8$, $m=6l$ and $9\le l$, or
$r=8$, $m=6l+1$ and $6\le l$. Using \ir{recall}
to restrict to those cases with $\alpha(I)=\beta(I)$
gives the result. 

\References

\bibitem{A} Alexander, J. {\it Singularit\'es 
imposables en position g\'en\'erale
aux hypersurfaces de \pr n}, Comp. Math. 68 (1988), 305--354.

\bibitem{AH1} Alexander, J. and Hirschowitz, A. {\it Polynomial
interpolation in several variables}, J. Alg. Geom. 4 (1995), 201--222.

\bibitem{AH2} \manyby. {\it Une lemme 
d'Horace diff\'erentielle:
application aux singularit\'es hyperquartiques de \pr 5}, J. Alg. Geom. 1 
(1992), 411--426.

\bibitem{AH3} \manyby. 
{\it La m\'ethode d'Horace 
\'eclat\'ee: application \`a l'interpolation 
en degr\'e quatre}, Inv. Math. 107
(1992), 585--602.

\bibitem{Bk} Baker, A. {\it A concise introduction to the theory
of numbers}, Cambridge University Press, 1984, xiii + 95 pp.

\bibitem{Bl} Ballico, E. {\it Generators for the homogeneous 
ideal of $s$ general points in \pr 3}, 
J.\ Alg.\  106 (1987), 46--52.

\bibitem{Cam} Campanella, G. {\it Standard bases of perfect homogeneous
polynomial ideals of height $2$}, 
J.\ Alg.\  101 (1986), 47--60.

\bibitem{Cat} Catalisano, M.\ V. {\it ``Fat'' points on a conic}, 
Comm.\ Alg.\  19(8) (1991), 2153--2168.

\bibitem{Ch} Chandler, K. {\it A brief proof of 
Alexander-Hirschowitz's Theorem},
preprint, 1996.

\bibitem{CM} Ciliberto, C. and Miranda, R. {\it On the dimension
of linear systems of plane curves with general multiple base points},
preprint, 1997.

\bibitem{DGM} Davis, E.\ D., Geramita, A.\ V., and 
Maroscia, P. {\it Perfect
Homogeneous Ideals: Dubreil's Theorems Revisited},
Bull.\ Sc.\ math., $2^e$ s\'erie, 108 (1984), 143--185.

\bibitem{EP}  Eisenbud, D. and Popescu, S. {\it Gale Duality 
and Free Resolutions of Ideals of Points}, preprint, 1996.

\bibitem{Fi} Fitchett, S. {\it Doctoral dissertation}, 
in preparation, 1996.

\bibitem{GO} Geramita, A. V. and Orrechia, F. {\it Minimally
generating ideals defining certain tangent cones}, J.\ Alg. 78
(1982), 36--57.

\bibitem{GM} Geramita, A. V. and Maroscia, P. {\it The ideal
of forms vanishing at a finite set of points of \pr n}, J.\ Alg. 90
(1984), 528--555.

\bibitem{GGR} Geramita, A.\ V., Gregory, D.\ and Roberts, L.
{\it Minimal ideals and points in projective space},
J.\ Pure and Appl.\ Alg. 40 (1986), 33--62.

\bibitem{Gi} Gimigliano, A. {\it Our thin knowledge of fat points},
Queen's papers in Pure and Applied Mathematics, no. 83,
The Curves Seminar at Queen's, vol. VI (1989).

\bibitem{\trans} Harbourne, B. {\it Complete linear 
systems on rational surfaces}, 
Trans.\ A.\ M.\ S.\ 289 (1985), 213--226. 

\bibitem{\vanc} \manyby. {\it The geometry of 
rational surfaces and Hilbert
functions of points in the plane},
Can.\ Math.\ Soc.\ Conf.\ Proc.\ 6 
(1986), 95--111.

\bibitem{\bowdoin} \manyby. {\it Automorphisms of K3-like
Rational Surfaces}, Proc.\ Symp.\ Pure Math.\ 46 
(1987), 17--28.

\bibitem{\ravello} \manyby. {\it 
Points in Good Position in \pr 2}, in:
Zero-dimensional schemes, Proceedings of the
International Conference held in Ravello, Italy, June 8--13, 1992,
De Gruyter, 1994.

\bibitem{\mtnwest} \manyby. {\it Rational 
surfaces with $K^2>0$}, Proc. A.M.S.
124 (1996), 727--733.

\bibitem{\antican} \manyby. {\it Anticanonical 
rational surfaces}, Trans. A.M.S. 349 (1997), 1191--1208.

\bibitem{\fatpts} \manyby. {\it Free Resolutions of Fat Point 
Ideals on \pr2}, preprint (available via my Web page), 1995, to appear,
JPAA.

\bibitem{Hi} Hirschowitz, A.
{\it Une conjecture pour la cohomologie 
des diviseurs sur les surfaces rationelles g\'en\'eriques},
Journ.\ Reine Angew.\ Math. 397
(1989), 208--213.

\bibitem{HS} Hirschowitz, A. and Simpson, C.
{\it La r\'esolution minimale de l'id\'eal d'un 
arrangement g\'en\'eral d'un
grand nombre de points dans \pr n}, Invent. Math. 126(1996), 467--503.

\bibitem{HSV} Hoa, L.T., Stuckrad, J., and Vogel, W. {\it Towards
a structure theory for projective varieties of degree $=$ 
codimension $+$ $2$}, 
J. Pure Appl. Algebra 71 (1991), 203--231.

\bibitem{I} Iarrobino, A. {\it Inverse system of a 
symbolic power, III: thin algebras
and fat points}, preprint 1994.


\bibitem{Lor} Lorenzini, A. {\it The minimal 
resolution conjecture}, J.\ Alg.
156 (1991), 5--35.

\bibitem{Ma} Manin, Y.\ I. Cubic Forms. North-Holland 
Mathematical Library
4, 1986.

\bibitem{Mu} Mumford, D. {\it Varieties defined by quadratic equations},
in: Questions on algebraic varieties, Corso C.I.M.E. 1969 Rome: Cremonese,
1969, 30--100.

\bibitem{N1} Nagata, M. {\it On rational surfaces, II}, 
Mem.\ Coll.\ Sci.\ 
Univ.\ Kyoto, Ser.\ A Math.\ 33 (1960), 271--293.

\bibitem{N2} \manyby. {\it On the 14-th problem of Hilbert}, 
Amer.\ J.\ Math.\ 33 (1959), 766--772.

\bibitem{O} Okuyama, H. {\it A note on conjectures 
of the ideal of $s$-generic points in \pr 4},
J. Math. Tokushima Univ. 25  (1991), 1--11.

\bibitem{Ra} Ramella, I. {\it An algorithmic 
approach to ideal generation of points},
Rend. Accad. Sci. Fis. Mat. Napoli (4) 56 (1989), 71--81.

\bye